\documentclass[11pt, reqno]{article}
\usepackage{amsfonts}
\usepackage{amsmath}
\usepackage{amsthm}
\usepackage{amssymb}
\usepackage{lastpage}
\usepackage{fancyhdr}
\usepackage{fancybox}
\usepackage{mathrsfs}

\def\b{\mathbb }
\def\phi{\varphi }

\swapnumbers

\theoremstyle{plain}
\newtheorem{theorem}{Theorem}[section]
\newtheorem{corollary}[theorem]{Corollary}
\newtheorem{lemma}[theorem]{Lemma}
\newtheorem{proposition}[theorem]{Proposition}
\theoremstyle{definition}
\newtheorem{definition}[theorem]{Definition}
\newtheorem{remark}[theorem]{Remark}

\numberwithin{equation}{section}

\begin{document}

\title{ Bessel convolutions on matrix cones:
         Algebraic properties and random walks}
\author{Michael Voit\\
michael.voit@math.uni-dortmund.de\\
Fachbereich Mathematik, Universit\"at Dortmund\\
          Vogelpothsweg 87\\
          44221 Dortmund, Germany}
\date{March 1, 2006}
\maketitle

\begin{abstract}
Bessel-type convolution algebras of bounded Borel
measures on the matrix cones of positive 
semidefinite  $q\times q$-matrices
over  $\b R, \b C, \b H$  were introduced recently by R\"osler.
These convolutions  depend on some continuous parameter,
 generate commutative hypergroup structures and have
 Bessel functions of matrix argument as characters.
Here, we first study the rich algebraic structure of
 these hypergroups. In particular, the subhypergroups and
 automorphisms are classified, and we  show
 that each quotient by a subhypergroup carries  a hypergroup
 structure of the same type.
  The algebraic properties  are
partially related to properties of random walks on 
 matrix Bessel hypergroups. In particular,  known properties of 
Wishart distributions, which form Gaussian convolution semigroups on these hypergroups,
are put into a new light. Moreover,
 limit theorems for random walks on these hypergroups are presented.
 In particular,  we obtain strong laws of large numbers and a central 
limit theorem with  Wishart distributions
as limits.
\end{abstract}


\section{Introduction} Recently, R\"osler
\cite{R} introduced positivity-preserving convolution algebras
 on the matrix cones $\Pi_q(\b F)$ of positive 
semidefinite  $q\times q$-matrices over  
$\b F= \b R, \b C, \b H$ which 
are related with Bessel functions 
of matrix argument and depend on  some continuous parameter.
With respect to this parameter, they interpolate the 
radial convolution algebras on non-squared matrix spaces
$M_{p,q}(\b F)$ with $p\geq q$ which are, to some extent, 
studied in  \cite {FT}. The convolutions of \cite{R} 
generate commutative hypergroup structures  on $\Pi_q(\b F)$
with Bessel functions of matrix argument as characters;
 see \cite{FT}, \cite{FK}
\cite{H} and \cite{Di} for matrix  Bessel functions and 
 \cite{BH}, \cite{J} for  hypergroups.
The present paper is devoted to algebraic and 
probabilistic aspects of these hypergroups.
In particular we shall show that these hypergroups
 admit  many subhypergroups and hypergroup
automorphisms. Moreover,  these algebraic properties  
 are closely related with probabilistic properties of random walks on 
these hypergroups. In particular,  
some known  properties of Wishart distributions will be seen under a new light.

 Before going into detail, we recall the  one-dimensional case:
For any dimension $p\ge1$, the  Banach-$*$-algebra $M_b(\b R^p)$ of 
bounded Borel measures on $\b R^p$  with the usual convolution
contains the space of all radial measures
$$M_b^{rad}(\b R^p):=\{ \mu\in M_b(\b R^p):\> u(\mu)=\mu
 \quad\text{ for all}\quad u\in O(p)\}$$
 as a Banach-$*$-subalgebra. If we identify
the set of all orbits under the standard action of the
 orthogonal group $O(p)$ on  $\b R^p$ 
with $\Pi_1:=[0,\infty[$ via  
$$p:x\mapsto |x|= (x_1^2 +\ldots +
  x_p^2)^{1/2},$$ 
then $p$ induces an isomorphism between the Banach spaces
  $M_b^{rad}(\b R^p)$ and  $M_b(\Pi_1)$. We thus may transfer the 
Banach-$*$-algebra  structure of  $M_b^{rad}(\b R^p)$ to  $M_b(\Pi_1)$
which inherits a commutative, associative, probability
preserving and weakly continuous convolution $*_p$.
Calculation in polar coordinates shows that for $p\ge 1$,
\begin{equation}\label{Bessel-faltung}
 \delta_r *_p \delta_s (f) = c_p\int_0^{\pi}
 f\bigl(\sqrt{r^2 +s^2 -2rs \cos\theta}\bigr) \sin^{p-2}\theta\, d\theta,
 \quad r,s \ge0, \, f\in C(\Pi_1)
\end{equation}
with a  normalization constant $c_p>0$ where for $p=1$  (\ref{Bessel-faltung})
degenerates to
\begin{equation}\label{cosinus-faltung}
\delta_r*\delta_s=\frac{1}{2}(\delta_{|r-s|}+\delta_{r+s})\quad\quad
(r,s\in\b R,\> r,s\ge0).
\end{equation}
The convolution on  $M_b(\Pi_1)$ is then obtained by linear,
weakly continuous extension. It is well-known that for all real $p>1$,
 Eq.(\ref{Bessel-faltung})
generates a commutative, probability
preserving, and weakly continuous convolution algebra on $M_b(\b R_+)$ which
interpolates
the integer cases. These convolutions
have no  group interpretation and are closely related with the 
known product formulas  
\[ j_\alpha(r) j_\alpha(s) = \delta_r *_p \delta_s(j_\alpha).\]
for  the  normalized Bessel functions
$\, j_\alpha(z) = \,_0F_1(\alpha +1; -z^2/4)\,$
with index $\alpha = p/2-1\ge -1/2$ (\cite{W}).
The space $\b R_+$  with the convolutions $*_p$ for $p\ge1$ provides a 
 prominent class of  commutative hypergroups, called 
Bessel-Kingman hypergroups (\cite{BH}).   
For $p>1$, these hypergroups have no nontrivial subhypergroups while in
the degenerated case 
$p=1$ the sets $c\b Z_+$ for  $c>0$ form the nontrivial subhypergroups.
Moreover, for $p\ge1$,  all hypergroup automorphisms are given by $x\mapsto cx$ for $c>0$,
see \cite{Z1}. Random walks on Bessel-Kingman hypergroups, i.e., 
Markov chains on $\b R_+$  with transition probabilities
given in terms of $*_p$  were investigated first by Kingman \cite{K}; for
the later development we refer to \cite{BH} and references therein. 
Kingman in particular obtained laws of large numbers and a central limit
theorem with a Rayleigh distribution as limit. For integers $p$, these
 limit theorems on $\b R_+$   are  just  
 radial reformulations of classical limit theorems on $\b R^p$.

We now turn to the higher rank case  in  \cite{R}.
For $p,q\in \b N$ with $p\geq q$, consider the space
$M_{p,q} = M_{p,q}(\b F)$  of $p\times q$ matrices over one of the division algebras
 $\b F = \b R, \b C$ or the quaternions $\b H$ with real dimension $d=1,2$ or $4$ respectively.
$M_{p,q}$ is a real Euclidean vector space of dimension $dpq$ with scalar product
$(x|y) = \mathfrak R \text{tr}(x^*y)$ where $x^* = \overline x^t$, 
$\mathfrak R t = \frac{1}{2}(t+ \overline t)$ is the real part of $t\in \b F$,  and
$\text{tr}$  the trace in $M_q(\b F):=M_{q,q}(\b F)$.
A  measure  on $M_{p,q}$ is called radial if it is 
invariant under the action of the unitary group $U_p= U_p(\b F)$ on $M_{p,q}$ 
by left multiplication,
\begin{equation}\label{leftaction}
U_p\times M_{p,q} \to M_{p,q}\,, \quad (u,x) \mapsto ux.
\end{equation}
This action is  orthogonal w.r.t.~the scalar product above, and   $x,y$
 are in the same $U_p$-orbit if and only if $x^*x = y^*y$.
Thus the  space of $U_p$-orbits is naturally parametrized by the 
cone $\Pi_q = \Pi_{q}(\b F)$ of positive semidefinite $q\times q$-matrices over $\b F$.
For $q=1$ and  $\b F= \b R$, we  have $\Pi_1=\b R_+$ and end up with
 the one-dimensional case above. We now use the projection 
$$ p: M_{p,q} \to \Pi_{q},\quad x\mapsto (x^*x)^{1/2},$$ 
with  the usual unique square root on  $\Pi_q$. 
Via this mapping
 the convolution algebra of radial measures on $M_{p,q}$ is transferred to 
 a commutative, associative, probability
preserving and weakly continuous convolution $*_p$ of measures on $\Pi_q$
which forms a commutative hypergroup. 
By construction (and results of \cite{FT}, \cite{H})  this
 convolution corresponds to a product formula for   Bessel functions $\mathcal J_\mu$ on the cone $\Pi_q$  with index $\mu=pd/2$. In \cite{R}, the convolution $*_p$ and  the product
 formula for the corresponding $\mathcal J_\mu$ is written down 
 in a way which allows for analytic continuation with respect to the index.
This leads to  ``interpolating'' commutative hypergroup structures $X_{q,\mu}$ on
$\Pi_q$ with a continuous real index $p\geq 2q$, i.e. $\mu
\geq d(q-1/2)$ and with matrix Bessel functions of index $\mu$ 
as hypergroup characters. These hypergroups are self-dual with the 
identity mapping as involution.
  The product formulas degenerate for $p=q, q+1\ldots,
2q$. For   non-integer $p\in ]q, 2q[$ there is unfortunately only a guess for
 explicit   product  formulas; see
\cite{R}. 

The present paper continues \cite{R}. In the first place, we
study algebraic properties of the matrix Bessel hypergroups of \cite{R}. 
 We first show that for each $a\in GL(q,\b F)$, the
map $T_a(r):=(ar^2a^*)^{1/2}$ is a hypergroup automorphism of $X_{q,\mu}$ 
 This
reveals that matrix Bessel hypergroups  in  higher rank admit a rich
structure of automorphisms similar to  the Euclidean spaces $\b F^d$.
The deepest result  will be the classification of all hypergroup
automorphisms $Aut(X_{q,\mu})$ for $\b F=\b R,\b C$. Indeed, we shall prove that for $\b
F=\b R$,  $Aut(X_{q,\mu})$ is the transformation group
$\{T_a:\> a\in GL(q,\b F)\}$, and that for $\b
F=\b C$ in addition the maps $\tau\circ T_a$ appear with $a\in GL(q,\b F)$
and $\tau$ the complex conjugation. We expect a similar result for $\b
F=\b H$, but we are unable to prove it here.
 In addition, we will classify all subhypergroups in the general case, that is 
for all $\b F$ and $\mu >d(q-1/2)$. More precisely, we
 prove that all subhypergroups are of   the form 
$$H_{k,u}:=\left\{
u\left(\begin{array}{cc}
\tilde r & 0\\0&0
\end{array}
\right)u^*: \> \tilde r\in\Pi_k\right\}$$
with $k\ge 0$ and a unitary matrix $u\in U_q$ (where $H_{0,u}=\{0\}$). 
We  also show that for fixed  $\mu$,  $H_{k,u}$ is 
canonically isomorphic with the hypergroup $X_{k,\mu}$, and that the
quotient  $X_{k,\mu}/H_{k,u}$ carries a  quotient hypergroup
structure and is isomorphic with $X_{q-k,\mu}$. The proofs of these algebraic
properties will rely more on the properties of matrix Bessel functions (which form
the hypergroup characters) rather than the explicit form of the convolution. 

The second part of this paper is devoted to probability theory on matrix
Bessel hypergroups. We there introduce convolution semigroups of probability
measures and  random walks and  show how Wishart
distributions fit into this concept. In particular, some known facts about 
 Wishart distributions and Wishart processes will appear under a new light.
Moreover, these Wishart distributions will appear as limits in a central limit
theorem. Besides this  central limit theorem, we  derive strong laws of
large numbers for random walks on $X_{q,\mu}$. The proofs of both types of
limit theorems are quite standard from a point of view of limit theorems on
commutative hypergroups and rely on the concept of  moment functions
which was developed by Zeuner \cite{Z1}, \cite{Z2} and others; see also the
monograph \cite{BH}. We point out that for the group case $\mu=pd/2$, all
limit theorems are just radial reformulations of classical limit theorems on 
the vector spaces $M_{p,q}$.

This paper  is organized as follows: In  Section 2 we collect some basic
  facts about matrix Bessel functions  and the corresponding
 hypergroups from  \cite{R}. Section 3  contains some general
facts about hypergroups which will be
useful for our concrete examples. In Section 4,  the
algebraic properties of matrix
Bessel hypergroups are studied. The remaining sections are then 
devoted to probability theory
on  matrix
Bessel hypergroups. 
Section 5 contains  basic properties of
convolution semigroups of probability measures, random walks,  and  Wishart
distributions. In Section 6 we then derive  a central limit
theorem as well as strong laws of large numbers for random walks on matrix cones.

\medskip\noindent
\emph{Acknowledgement.}  It is a pleasure to thank Margit R\"osler for many
discussions and comments.

\section{Bessel convolutions on matrix cones}

Here we collect some basic
notions and  facts about matrix Bessel functions
 and  matrix
Bessel hypergroups from \cite{FK}, \cite{BH}, \cite{J}, \cite{R}.

\subsection{Bessel functions associated with matrix cones}

Let $\b F$ be one of the real division algebras $\b F= \b R, \b C$ 
or  $\b H$ with real dimension $d=1,2$ or $4$ respectively.
Denote the usual conjugation in $\b F$ by $t\mapsto \overline t$,
the real part  of $t\in \b F$ by $\mathfrak R t = \frac{1}{2}(t + \overline
t)$, and $|t| = (t\overline t)^{1/2}$ its norm.

For $p,q\in \b N$ we denote  by $M_{p,q}:=M_{p,q}(\b F)$ 
the vector spaces of all  $p\times q$-matrices
 over $\b F$, and we put $M_{q}:=M_{q,q}$. 
We consider the set 
\[H_q = H_q(\b F)= \{x\in M_{q}(\b F): x= x^*\}\]
of Hermitian $q\times q$-matrices over $\b F$ as a Euclidean
 vector space with scalar product 
 $(x\vert y) := \mathfrak R {tr}(x^*y)$  and 
the associated norm $\|x\|=(x\vert x)^{1/2}$. Here 
$x^* := \overline x^t$ and ${tr}$ denote the trace.
Its dimension is given by  $dim_{\b R}H_q=n:= q + \frac{d}{2}q(q-1)$.
Let further
\[\Pi_q:=\{x^2:\> x\in H_q\} = \{ x^*x: x\in H_q\}\]
be the set of all positive semidefinite
matrices in $H_q$, and $\Omega_q$ its topological 
interior which consists of all strictly positive
definite matrices. $\Omega_q$ is a symmetric cone, i.e. an open convex
self-dual cone which  whose  linear automorphism group acts transitively, see \cite{FK} for details.

To define the Bessel functions associated with the symmetric cone $\Omega_q$ we first  introduce their basic building blocks,
the so-called spherical polynomials. These are just 
 the polynomial spherical functions of $\Omega_q$ considered as a Riemannian symmetric space. They are indexed
by partitions $\lambda = (\lambda_1 \geq \lambda_2\geq \ldots \geq \lambda_q)
\in \b N_0^q$ (we write $\lambda \geq 0$ for short) and given by 
\[ \Phi_\lambda (x) = \int_{U_q} \Delta_\lambda(uxu^{-1})du, \quad x\in H_q\]
where $du$ is the normalized Haar measure of $U_q$ and $\Delta_\lambda$ is the power function on  $H_q$ with
\[ \Delta_\lambda(x) := \Delta_1(x)^{\lambda_1-\lambda_2} \Delta_2(x)^{\lambda_2-\lambda_3} \cdot\ldots\cdot \Delta_q(x)^{\lambda_q};\]
the $\Delta_i(x)$ are the principal minors of the determinant $\Delta(x)$, see
\cite{FK} for details. The $\Phi_\lambda$ are homogeneous of degree 
$|\lambda|=\lambda_1 +\ldots + \lambda_q$. 
There is a renormalization
$Z_\lambda = c_\lambda \Phi_\lambda$
with constants $c_\lambda >0$ depending on the underlying cone
such that
\begin{equation}\label{traceid}
(tr \,x)^k \,=\, \sum_{|\lambda|=k} Z_\lambda(x)
\quad\quad{\rm for}\>\> k\in \b N_0,
\end{equation}
see Section XI.5.~of \cite {FK} (the $Z_\lambda$ are called
zonal polynomials there). By construction, they are invariant under conjugation by $U_q$ and thus depend only on the eigenvalues of their argument.
More precisely, for $x\in H_q$ with eigenvalues $\xi = (\xi_1, \ldots, \xi_q)\in \b R^q$,
\[ Z_\lambda(x) = C_\lambda^\alpha(\xi), \quad \alpha = \frac{2}{d}\]
where the $C_\lambda^\alpha$ are the Jack polynomials of index
$\alpha$ in a suitable normalization (c.f. \cite {FK}, \cite{Ka}, \cite{R}).
They are homogeneous of degree $|\lambda|$ and symmetric in their arguments.

The matrix Bessel functions associated with the cone $\Omega_q$ are now defined
as  $_0F_1$-hyper\-geometric series in terms of the $Z_\lambda$, as follows: 
\[ \mathcal J_\mu(x) =
 \sum_{\lambda\geq 0} \frac{(-1)^{|\lambda|}}{(\mu)_\lambda|\lambda|!} Z_\lambda(x), \quad x\in H_q \]
where for $\lambda=(\lambda_1, \ldots \lambda_q)\in \b N_0^q$,
 the generalized Pochhammer symbol $(c)_\lambda$ is 
\[ (\mu)_\lambda =\, (\mu)_\lambda ^{2/d}\,=\,\prod_{j=1}^q \bigl(\mu-\frac{d}{2}(j-1)\bigr)_{\lambda_j}\]
and 
$\mu\in \b C$ is an index with $(\mu)_\lambda \not= 0$ for  $\lambda \geq 0.$ 
This  series converges absolutely for  $x\in H_q$, see \cite{FK}. Later on, we
need the linear terms in the expansion of $\mathcal J_\mu$. By \eqref{traceid} we have
\begin{equation}\label{homexp}
 \mathcal J_{\mu}(x) = 1 - \frac{1}{\mu}\, tr\, x + O(\|x\|^2).\end{equation}

To describe the main results of \cite{R}, we need the notion of
a hypergroup, which will be introduced in the following section.

\subsection{Hypergroups}
A hypergroup $(X,*)$
consists of a locally compact Hausdorff space $X$ and a multiplication $*$,
called  convolution, on the Banach space $M_b(X)$ of all bounded
regular complex Borel measures  with the total variation norm as norm,
such that $(M_b(X),*)$ becomes a Banach algebra, and such that  $*$ is
weakly continuous and probability preserving and  preserves compact
supports of measures.
 Moreover, there  exists an identity $e\in X$
with  $\delta_e*\delta_x=
\delta_x*\delta_e=\delta_x$ for $x\in X$, as well as
  a continuous involution $x\mapsto\bar x$ on $X$ such that for $x, y\in X$,
$e\in supp(\delta_x*\delta_y)$ is equivalent to $ x=\bar y$,
and
$\delta_{\bar x} * \delta_{\bar y} =(\delta_y*\delta_x)^-$.
Here for $\mu\in M_b(X)$, the measure  $\mu^-$ is given by
$\mu^-(A)=\mu(A^-)$ for Borel sets $A\subset X$.

A hypergroup $(X,*)$ is called commutative if and only if so is the
convolution $*$, and hermitian, if the hypergroup involution is the identity.
 Hermitian hypergroups are  commutative.

It is well-known that each commutative hypergroup admits 
 a (up to normalization)  unique
Haar measure $\omega\in M^+(X)$ which is characterized by
$\, \omega(f)=\int_X f(x*y)\> d\omega(y)$ for all $x\in X$ and 
all compactly supported,
 continuous functions  $\,f\in C_c(X)$
where we use the notation
\[f(x*y):= \int_X f\> d(\delta_x*\delta_y).\]
Similar to the dual  of a locally compact abelian group,
we define the dual 
$$\widehat X:=\{\alpha\in C_b(X):\>\alpha\ne 0,\>\>  \alpha(\bar x*y)=
 \overline{\alpha(x)}\alpha(y)\>\>\text{for all}\>\> x,y\in X\}.$$
$\widehat X$ is a locally compact Hausdorff space w.r.t.
 the topology of
  compact-uniform convergence, and  its elements are called characters.
The Fourier transform on $ L^1(X, \omega)$ is defined by
 $\widehat f(\alpha):= \int_X f(x)\overline{\alpha(x)}\>
d\omega(x)$, $\alpha\in \widehat X$, and the Fourier-Stieltjes transform of
measures is defined in the same way. It is well-known that for a fixed Haar
 measure $\omega$ on $X$ there is a unique Plancherel measure $\pi$ on $\hat
 X$ such that the Fourier transform becomes an $L^2$-isometry between
 $L^2(X,\omega)$ and $L^2(\hat X,\pi)$.

Interesting examples of  hypergroups are given as follows. 

\subsection{Bessel convolutions on matrix cones}
For natural numbers $p, q$, consider  the matrix space $M_{p,q}=M_{p,q}(\b F)$ of 
$p\times q$-matrices over $\b F$.
 We regard $M_{p,q}$ as a real vector space with the Euclidean scalar product
 $(x\vert y) := \mathfrak R\text{tr}(x^*y)$ and norm $\,\|x\| =
 \sqrt{\text{tr}(x^*x)}$.
 Consider the action of the
 unitary group $U_p$  on $M_{p,q}$ by left multiplication,  
\[ U_p\times M_{p,q} \to M_{p,q}\,, \quad (u,x) \mapsto ux.\] 
This action is  orthogonal w.r.t.~the scalar product above, and   $x,y$
 are in the same $U_p$-orbit if and only if $x^*x = y^*y$. The space
$M_{p,q}^{U_p}$ of all orbits for this action
can therefore be identified with the space 
$\Pi_q= \Pi_q(\b F)$ of positive semidefinite $q\times q$ matrices over $\b F$ via 
\[ U_p.x \mapsto \,\sqrt{x^*x}\, =:\, |x|.\]
 Here for $r\in \Pi_q\,,\, \sqrt{r}$ is the unique positive semidefinite
 square root of $r$.  
This  bijection is a homeomorphism  w.r.t. the quotient topology on
 $M_{p,q}^{U_p}$.

Now consider the Banach-$*$-algebra
$$M_b^{rad}(M_{p,q}):=\{\mu\in M_b(M_{p,q}):\> u(\mu)=\mu \quad \text{for all}\>\> u\in U_p\}$$
of all radial regular Borel   measures  on $M_{p,q}$, and the canonical
projection 
$$ p: M_{p,q} \to \Pi_{q},\quad x\mapsto (x^*x)^{1/2},$$ 
with  the usual unique square root on  $\Pi_q$. 
Via this mapping, the convolution on $M_b^{rad}(M_{p,q})$  is transferred to 
 a commutative, associative, probability
preserving and weakly continuous convolution $*_p$ of measures on $\Pi_q$
which forms a commutative hypergroup. 
By construction (and results of \cite{FT}, \cite{H})  this
 convolution corresponds to a product formula for Bessel functions $\mathcal
 J_\mu$ on  $\Pi_q$  with index $\mu=pd/2$. In \cite{R}, the convolution $*_p$ and  the product
 formula for the corresponding $\mathcal J_\mu$ is written down explicitly
 (see Eq.~(\ref{def-convo}) below)
 in a way which allows for analytic continuation with respect to the index $\mu$.
This leads to  ``interpolating'' commutative hypergroup structures $X_{q,\mu}$
on $\Pi_q$ for all indices $\mu> \rho-1$ with
\[ \rho:= d\bigl(q-\frac{1}{2}\bigr) +1.\]
For $\mu\le\rho-1$ having the form $\mu=pd/2$ with $p\in\b N$, there exist
also degenerated versions of the product formula (\ref{def-convo}) below; it
is however not clear at the moment whether these discrete cases can be
embedded into a continuous families  of convolution and product formulas; see
the discussion in \cite{R}. 

In the following we use the  abbreviations 
$ D_{q} = \{ v\in M_{q,q}: v^*v < I\}$
and 
$$ \kappa_\mu := \int_{D_q} \Delta(I-v^*v)^{\mu-\rho} dv.$$
of \cite{R}, where an explicit formula for $\kappa_\mu$ is given in \cite{R}.
The following result contains some of the main results of  \cite{R}.

\begin{theorem}\label{main2}
Let $\mu\in\b R$ with $\mu > \rho-1.$ Then
\parskip=-1pt
\begin{enumerate}\itemsep=-1pt
\item[\rm{(a)}] The assignment
\begin{equation}\label{def-convo}
(\delta_r *_\mu \delta_s)(f) := 
\frac{1}{\kappa_\mu}\int_{D} f\bigl(\sqrt{r^2 + s^2 + svr + rv^*\!s}\,\bigr)\,
\Delta(I-vv^*)^{\mu-\rho}\, dv
\end{equation}
for $f\in C_c(\Pi_q)$ defines a hermitian hypergroup structure on $\Pi_q$ with 
neutral element $0\in \Pi_q$. 
The support of $\delta_r*_\mu\delta_s$ satisfies
\[\text{supp}(\delta_r*_\mu\delta_s) \subseteq  
\{t\in \Pi_q: \|t\|\leq \|r\|+\|s\|\}.\]
\item[\rm{(b)}] A Haar measure of this hypergroup 
$X_{q,\mu}:= (\Pi_q, *_\mu)$ is given by
\[ \omega_\mu(f) = 
\frac{\pi^{q\mu}}{\Gamma_{\Omega_q}(\mu)}
\int_{\Omega_q} f(\sqrt{r}) \Delta(r)^\gamma dr\]
with $\,\gamma = \mu-\frac{d}{2}(q-1)-1\,=\, \mu-\frac{n}{q}.$
\item[\rm{(c)}]
The dual space of $X_{q,\mu}$ is given by
\[\widehat{X_{q,\mu}} = \,\{\phi_s = \phi_s^\mu : s\in \Pi_q\},\]
 and $X_{q,\mu}$ is self-dual via the homeomorphism $s\mapsto \phi_s$.
Under this identification of $\widehat{X_{q,\mu}}$ with $X_{q,\mu}\,$,
  the Plancherel measure on $X_{q,\mu}\,$ is $(2\pi)^{-2\mu q}\omega_\mu$.
\end{enumerate}
\end{theorem}

Notice that in our normalization of the Haar  measure  for
$\mu=pd/2$, $ \omega_\mu$ is just the image of the
 Lebesgue measure on $M_{p,q}\simeq\b R^{dpq}$ under the canonical projection
 $p:M_{p,q}\simeq\b R^{dpq}\to\Pi_q$.

The explicit  product formula (\ref{def-convo}) may be used to describe the
support of convolution products in more detail. Some special  cases will be 
considered in Lemmas \ref{support} and \ref{step6} below.

\section{Some results about hypergroups}

In this section we  collect  some further general notions and results about commutative hypergroups. 
We also prove some more or less straightforward results about  hypergroup automorphisms which will be
useful  for the Bessel convolutions below. 

We start with the following simple observation will be useful  several times.

\begin{lemma}\label{hypergroup-property1}
If a locally compact space $X$ carries two commutative
  hypergroup convolutions $*_1$ and  $*_2$ such that the dual $(X,*_1)^\wedge$
  is contained in $(X,*_2)^\wedge$, then $*_1=*_2$.
\end{lemma}

\begin{proof}
 Let $\alpha\in(X,*_1)^\wedge\subset(X,*_2)^\wedge$. Then, for
  all $x,y\in X$,
$$(\delta_x*_1 \delta_y)^\wedge(\alpha)= \bar\alpha(x) \bar\alpha(y)=
 (\delta_x*_2 \delta_y)^\wedge(\alpha) .$$
As the Fourier-Stieltjes transform of measures on $(X,*_1)$ is injective (see 
\cite{J} or \cite{BH}), we obtain $\delta_x*_1 \delta_y=\delta_x*_2 \delta_y$ for  $x,y\in X$.
\end{proof}

\begin{definition}
Let $X,Y$ be commutative hypergroups. A closed set $H\subset X$ is called a
subhypergroup, if for all $x,y\in H$, we have $\bar x\in H$ and
$\{x\}*\{y\}:=supp(\delta_x*\delta_y)\subset H$.
Moreover, a continuous 
mapping $T:X\to Y$ is  a hypergroup homomorphism, if $T(\delta_x* \delta_{\bar
  y})=
\delta_{T(x)}* \delta_{\overline{T(y)}}$  for all $x,y\in X$ where the mapping
$T:X\to Y$ is extended to bounded Borel measures by taking images of measures.
 The notions of a hypergroup isomorphisms and automorphisms are similar.
\end{definition}

We need the following observations:

\begin{lemma}\label{hypergroup-property2}
Let $(T_k)_{k\ge 1}$ be a sequence of hypergroup automorphisms on a
  hypergroup $(X,*)$ which converges pointwise to some continuous mapping
  $T:X\to X$ such that $T(X)$ is closed in $X$. Then $T$ is a hypergroup homomorphism
  from $X$ onto the subhypergroup $T(X)$.
\end{lemma}

\begin{proof}
 Let $x,y\in X$ and $f\in C_b(X)$. Then, by the weak continuity
  of the convolution,  
$$\lim_k \int f\> dT_k(\delta_x* \delta_{\bar y}) = \lim_k \int f\>
 d(\delta_{T_k(x)}*\delta_{\overline{T_k(y)}})=
\int f\> d(\delta_{T(x)}*\delta_{\overline{T(y)}}).$$
On the other hand, by the transformation formula and dominated convergence,
$$\lim_k \int f\> dT_k(\delta_x* \delta_{\bar y}) = \lim_k \int f\circ T_k\>
 d(\delta_{x}*\delta_{\bar y})=\int f\circ T\>d(\delta_{x}*\delta_{\bar y})
=\int f\> dT(\delta_{x}*\delta_{\bar y}).$$
This proves $T(\delta_x* \delta_{\bar y})=\delta_{T(x)}*
 \delta_{\overline{T(y)}}$, and that $T(X)$ is a subhypergroup.
\end{proof}

\begin{proposition}\label{prop-automor}
Let $X$ be a commutative hypergroup and  $T\in
Aut(X)$. Then:
\begin{enumerate}\itemsep=-2pt
\item[\rm{(1)}] The Haar measure   $\omega$  of $X$ satisfies
  $T(\omega)=c_T\omega$ for some constant $c_T>0$.
\item[\rm{(2)}] For $f,g\in C_c(X)$ and $x\in X$, $((f\circ T)*(g \circ T))(x)=
c_T\cdot (f*g)(T(x))$.
\item[\rm{(3)}] There exists a dual homeomorphism $T^*:\hat X\to\hat X$ with
$(T^*(\alpha))(x)=\alpha(T(x))$ for $\alpha\in\hat X, x\in X$. This mapping
$T^*$ maps the support $supp\> \pi\subset\hat X$ of the Plancherel measure of
$X$
onto  $supp\> \pi$. 
\item[\rm{(4)}] If $\hat X$ carries a dual hypergroup structure $(\hat X, *)$ which
is uniquely determined by
$$\int_{\hat X}\gamma(x)\> d(\delta_\alpha*\delta_\beta)(\gamma)=
\alpha(x)\beta(x)
\quad{\text for}\quad \alpha,\beta\in\hat X,\> x\in X,$$
then $T^*$ as defined in (2) is a hypergroup automorphism on $\hat X$.
\end{enumerate}
\end{proposition}

\begin{proof}
\begin{enumerate}\itemsep=-2pt
\item[\rm{(1)}] For $f\in C_c(X)$ and $x\in X$, 
\begin{align}T(\omega)(f_x):&= \int_X
  f(x*y) \> dT(\omega)(y)=\int_X f(x*T(w)) \> d\omega(w)\nonumber\\&=
\int_X f\circ T(T^{-1}(x)*w) \> d\omega(w)=\int_X f\circ T(w)\>
  d\omega(w)=T(\omega)(f).\nonumber
\end{align}
Thus $T(\omega)$ is  a Haar measure. The uniqueness of $\omega$ (\cite{J}) now yields the claim.
\item[\rm{(2)}] For $x\in X$,
\begin{align}((f\circ T)*(g \circ T))(x)&=\int_X (f\circ T)(x*\bar y)(g \circ
  T)(y)\> d\omega(y)
\nonumber\\ =& \int_X f(T(x)*\bar{T(y)}) g(T(y))\> d\omega(y)\nonumber\\
&=\int_X  f(T(x)* \bar w)g(w) \> dT(\omega)(w)=c_T\cdot (f*g)(T(x)).\nonumber
\end{align}
\item[\rm{(3)}] As for $x,y\in X$ and  $\alpha\in\hat X$
\begin{align}\int_ \alpha(T(z))\> d(\delta_x*\delta_{\bar y})(z)&=
\int_X \alpha(w) \> d(T(\delta_x*\delta_{\bar y}))(w)
 \nonumber\\ & =\int_X \alpha(w) \> d(\delta_{T(x)}*\delta_{\bar T(y)}))(w)=
\alpha(T(x))\overline{\alpha(T(y))},\nonumber
\end{align}
 $(T^*(\alpha))(x)=\alpha(T(x))$   defines for for  $\alpha\in\hat X$
a character $T^*(\alpha)\in \hat X$. Moreover, $T^*:\hat X\to \hat X$ is
obviously a homeomorphism (remember that $\hat X$ is equipped with the
topology of locally uniform convergence).

In order to check $T^*(supp\>\pi)\subset supp\>\pi$, take $\alpha\in
supp\>\pi$. Then, by \cite{V0}, there is a sequence $(f_n)_{n\ge1}\subset
C_c(X)$ with $f_n*f_n^*\to\alpha$ locally uniformly. Thus,
$$T(c_T^{-1/2}\cdot f_n)* T(c_T^{-1/2}\cdot f_n)^* = T(f_n*f_n^*)\to
T(\alpha)$$
locally uniformly, which conversely implies $T(\alpha)\in supp\>\pi$ by \cite{V0}.
\item[\rm{(4)}] For all $x\in X$ and $\alpha,\beta\in \hat X$ we have
$$T^*(\alpha*\beta^*)(x)= \int_{\hat X} \gamma(T(x))\> d(\delta_\alpha
*(\delta_{\bar\beta})(\gamma)=
\alpha(T(x))\overline{\beta(T(x))}$$
and
\begin{align} 
(T^*(\alpha)*\overline{T^*(\beta)})(x)&= \int_{\hat X} \gamma(x)\>
d(\delta_{T^*(\alpha)} *(\delta_{\overline{T^*(\beta)}})(\gamma)=
 T^*(\alpha)(x)\cdot
 \overline{T^*(\beta)(x)}\nonumber
\\\nonumber &=\alpha(T(x))\overline{\beta(T(x))}.\nonumber
\end{align}
\end{enumerate}
\end{proof}

\section{Automorphisms, subhypergroups, and quotients}

In this section we collect algebraic properties of the Bessel convolutions on
matrix cones. We  fix parameters $d,q,\mu,\rho$ as in Section in Section 2 and consider the associated hypergroup
structure $X_{q,\mu}$ on $\Pi_q$. 
One major task will be the classification of all hypergroup automorphisms of $X_{q,\mu}$.
For this, we first determine a group of hypergroup
automorphisms.
For this, we recall  that $GL(q):=GL(q, \b F)$
acts on $\Pi_q$ as a group of homeomorphisms via
$$T_a(r):=\sqrt{ar^2a^*} \quad\quad{\rm for}\quad a\in GL(q),\> r\in\Pi_q.$$
To check that the $T_a$ are in fact hypergroup  automorphisms, we  observe:

\begin{lemma}\label{trafo-charakter}
Let $s,r\in\Pi_q$ and $a\in M_q$. Then $\phi_s(T_a(r))=\phi_{T_{a^*}(s)}(r)$.
\end{lemma}

\begin{proof}
As the Bessel function ${\mathcal J}_\mu(r)$ depends on the spectrum  of $r\in\Pi_q$
only, we have
$$\phi_s(T_a(r))= {\mathcal J}_\mu(\frac{1}{4} sar^2a^*s)=
 {\mathcal J}_\mu(\frac{1}{4}a^*s^2ar^2)=\phi_{T_{a^*}(s)}(r).$$
\end{proof}

\begin{corollary}\label{Fourier-auto}
Let  $\mu\in M_b(\Pi_q)$,  $a\in M_q$, and  $s\in \Pi_q$. Then
 $\widehat{T_a(\mu)}(s)=\hat\mu(T_{a^*}(s))$.
\end{corollary}

\begin{proof}
$\widehat{T_a(\mu)}(s)=\int \phi_s(r)\> dT_a(\mu)(r)=
\int\phi_s(T_a(r))\>d\mu(r)= \int\phi_{T_{a^*}(s)}(r)\>d\mu(r)
=\hat\mu(T_{a^*}(s)).$
\end{proof}

\begin{proposition}\label{check-auto}
$\{T_a:\> a\in GL(q)\}$ is a group of hypergroup automorphisms of   $X_{q,\mu}$.
\end{proposition}

\begin{proof}
Fix $a\in GL(q)$. Then the homeomorphism $T_a$ induces a further commutative
hypergroup structure $(\Pi_q,*_a)$ on $\Pi_q$ by 
$$\delta_x*_a \delta_y:=T_a(\delta_{T_a^{-1}(x)}* \delta_{T_a^{-1}(y)})
\quad\quad (x,y\in \Pi_q.$$ 
It is easy to check that 
the dual
space of this  commutative hypergroup is 
 $\{\phi_s\circ T_{a^{-1}}:\> s\in \Pi_q\}$ where this space agrees
 with the dual
of $X_{q,\mu}$ by Lemma \ref{trafo-charakter}.  Lemma
\ref{hypergroup-property1} now  shows that the hypergroups $(\Pi_q,*_a)$ and 
$X_{q,\mu}$ agree, and hence $T_a$ is a hypergroup automorphism. As the $T_a$
obviously form a group, the proof is complete.
\end{proof}

 The preceding proposition may be also checked directly via the
explicit product formula (\ref{def-convo}), but in our eyes this approach is
more involved. Our approach has the further advantage that it works also for 
$\mu\le\rho-1$ in which  Eq.~(\ref{def-convo}) has a degenerated form (or is even unknown).

We prove below that for $\b F= \b R$  and 
$\mu>\rho-1$, the group $\{T_a:\> a\in GL(q)\}$
is  the  group of all hypergroup automorphisms of $X_{q,\mu}$.
This is however not correct for  $\b F=\b C, \b H$ and $q\ge 2$.
For instance, for $\b F=\b C$ and $q\ge2$, complex conjugation on $\Pi_q$ is an
automorphism which is not of the form above; for details see Theorems \ref{classiauto-reell}
 and \ref{classiauto-complex}  below.

We next determine all subhypergroups. For this we need:

\begin{lemma}\label{support}
Let $\mu>\rho-1$. Then for all $r\in\Pi_q$ and $c\in]0,1]$,
$$\{r\}*\{cr\}:=supp(\delta_r*\delta_{cr})=\{s\in\Pi_q:\> 
(1-c)r\le s\le (1+c)r\}.$$
\end{lemma}

\begin{proof}
We find  a suitable automorphism $T_a$ which maps $r$ into the  diagonal matrix
$I_j:=diag(1,\ldots,1,0,\ldots,0)$ with rank $j:=rank\> r$. We therefore
 may assume without loss of generality  $r=I_j$. 

By Eq.~(\ref{def-convo}), we have 
$$\{I_j\}*\{cI_j\}=\{
  \sqrt{(1+c^2)I_j+c(v+v^*)}\in\Pi_q:\> v\in M_q, \> vv^*\le I_j\} .$$ 
To simplify this set, we observe 
for $s\in\Pi_q$ and $h\in H_q$ that
$$s=\sqrt{(1+c^2)I_j+2ch} \quad\Longleftrightarrow \quad
h=(s^2-((1+c^2)I_j)/(2c),$$
where for $s,h$ coupled in this way, $(1-c)I_j\le s\le (1+c)I_j$
 is equivalent to
$-I_j\le h\le I_j$. Therefore, for a given $s$ with this property
we may take $v=h\in H_q$ and obtain
the inclusion $\supset$ in the statement of the lemma.
Conversely, for $v\in M_q$ with $vv^*\le I_j$, the
spectral norm of $v$ is bounded by 1; hence the spectral norm of
$h:=(v+v^*)/2$ is also  bounded by 1 which means for this hermitian matrix 
$-I_j\le h\le I_j$. This proves the converse inclusion 
 and completes the proof.
\end{proof}

\begin{remark} For $\mu=\rho-1$, the statement of Lemma \ref{support} is no
  longer correct. In fact, the degenerated explicit product formula 
in Proposition 3.16  of \cite{R} and 
some matrix computation shows
 that here for instance, for the identity matrix $I$ the set
 $\{I\}*\{I\}$ consists of those $s\in\Pi_q$ with eigenvalues 
$\lambda_1,\ldots,\lambda_q\in [0,2]$ with $\sum_{i=1}^q
 (\lambda_i^2/2-1)^2\ge1$, i.e., $\{I\}*\{I\}$ contains a hole.

In particular, for $q=1$ and $\b F=\b R$ (i.e. $d=1$) and $\mu=\rho-1=1/2$ we
have the  degenerated Bessel convolution (\ref{cosinus-faltung}).
This cosine hypergroup on $[0,\infty[$ has the discrete subhypergroups
$c\b N_0$ for 
$c>0$. This example shows in particular 
that in the following proposition we partially must
restrict our attention to the case $\mu>\rho-1$. 
\end{remark}

\begin{proposition}\label{all-subhyper}
Let $\mu\ge \rho-1$, $k\in\{0,1,\ldots,q\}$, and $u\in U_q$. Then
$$H_{k,u}:=\left\{
u\left(\begin{array}{cc}
\tilde r&0\\0&0
\end{array}
\right)u^*: \> \tilde r\in\Pi_k\right\}$$
(with $H_{0,u}=\{0\}$) is a subhypergroup of $X_{q,\mu}$, and the mapping
$\tilde r\mapsto u\left(\begin{array}{cc}
\tilde r&0\\0&0
\end{array}
\right)u^*$
is a hypergroup isomorphism between $X_{k,\mu}$ and $H_{k,u}$.
Moreover, for $\mu > \rho-1$ all subhypergroups of $X_{k,\mu}$ are given in this way.
\end{proposition}

\begin{proof}
The $H_{k,I}$ are obviously subhypergroups by 
 Eq.~(\ref{def-convo}). Using suitable
automorphisms, it becomes also clear that the  $H_{k,u}$ are subhypergroups
for arbitrary $u$.
In order to check the claimed isomorphism, we may  assume $u=I$ without
loss of generality.
It is possible to check the isomorphism property directly via Eq.~(\ref{def-convo})
and some computation. We here however prefer 
 another method and  observe 
 that the Jack polynomials  $C_\lambda^\alpha$ in $q$ and $k$ variables
 respectively
satisfy
$C_\lambda^\alpha(0,\ldots,0,\xi_1,\ldots,\xi_k)=C_\lambda^\alpha(\xi_1,\ldots,\xi_k)$
for $0\le\xi_1\le \ldots\le\xi_k$ by their very definition; 
see Stanley \cite{S}
Hence, by the definition of the $q$- and $k$-dimensional Bessel functions respectively, 
 \begin{equation}{\mathcal J}_\mu^q\left(\left(\begin{array}{cc}
\tilde r&0\\0&0
\end{array}
\right)\right)= {\mathcal J}_\mu^k(\tilde r)\quad\quad {\rm for}\quad\quad \tilde r\in
\Pi_k.
\end{equation}
Therefore,  all characters of $X_{k,\mu}$ appear as restrictions of
characters on $X_{q,p}$ to $H_{k,u}$ where these restrictions
are obviously characters on $H_{k,u}$. Lemma \ref{hypergroup-property1} now
shows that the hypergroup structures on $H_{k,u}$ and $X_{u,\mu}$ are
equal as claimed.

We still have to show that for $\mu > \rho-1$, all subhypergroups appear as
some $H_{k,u}$. For this we show that each  subhypergroup of
$X_{q,\mu}$ which is not contained in some $H_{q-1,u}$ with $u\in U_q$ must be
equal to $X_{q,\mu}$. As this says that each proper subhypergroup is contained 
in some $H_{q-1,u}$, we conclude from the first part of the proposition and 
induction that each subhypergroup appears as some $H_{k,u}$.
In order to prove the claim above, consider some 
subhypergroup $H$  which is not contained in some $H_{q-1,u}$. Let $a\in H$ be
an arbitrary element with rank  $k<q$ where we may assume without loss of
generality $a\in H_{k,I}$ after using a suitable automorphism. We then find
some
$b=\left(\begin{array}{cc}
 *&*\\ * & c
\end{array}
\right)\in H$ with $c\in \Pi_{q-k}$ and $c\ne0$. Then, by 
Eq.~(\ref{def-convo}), $a+b\in \{a\}*\{b\}\subset H$, and $a+b$ has rank at
least $k+1$. Iterating this argument, we find some $r\in H$ with full rank. 
Lemma \ref{support} now shows   that $\{r\}*\{r\}\subset H$ contains a
neighborhood $U$ of $0$ in $\Pi_q$. Applying  \ref{support} to elements of $U$
several times, finally implies $H=\Pi_q$ as claimed.
\end{proof}

\begin{remark}\label{quotient}
Let $a\in M_q$ be a matrix with rank $k\in\{0,\ldots,q\}$. 
We find $u,v\in U_q$ and a diagonal matrix $b=diag(b_1,\ldots,b_k,0,\ldots,0)$
with $b_1,\ldots,b_k\ne0$ such that $a=ubv^*$ holds. The mapping 
$T_a$  with $T_a(r):=\sqrt{ar^2a^*}$  then obviously is a
continuous and open mapping from $\Pi_q$ onto the  subhypergroup $H_{k,u}$.
Moreover, $T_a$ is a hypergroup homomorphism. To check this, 
choose a sequence
$(a_k)_k\subset GL(q)$ with $a_k\to a$. As then $T_{a_k}\to T_a$ pointwise on
$X_{q,\mu}$, the assertion follows from  Proposition \ref{check-auto}
 and Lemma \ref{hypergroup-property2}. We next notice that the kernel of
 $T_a$ is 
$$kern\> T_a :=\{r\in\Pi_q:\> T_a(r)=0\}\> =\> 
\left\{v\left(\begin{array}{cc}
 0&0\\0&\tilde r
\end{array}
\right)v^*:\> \tilde r\in \Pi_{q-k}\right\}
\> =\> H_{q-k,\tilde v}$$
for  
$$\tilde v:=\left(\begin{array}{ccccc}
0&0&\cdots&0&1 \\ 0&0&\cdots&1&0 \\ &&\cdots&&\\ 1&0&\cdots &0&0
\end{array}\right)\cdot v \in U_q,$$ and that $kern \> T_a$ is a subhypergroup
isomorphic with $X_{q-k,\mu}$ by Proposition \ref{all-subhyper}.

Now let $H_{k,u}$ ($0\le k\le q$, $u\in U_q$) be an arbitrary subhypergroup of
$X_{q,\mu}$ As $H_{k,u}$ appears as kernel of some hypergroup homomorphism $T$
from $X_{q,\mu}$ onto $X_{q-k,\mu}$ by the preceding considerations, we
conclude from  abstract results on hypergroup homomorphisms (see for
instance \cite{V1}) that the quotient space 
$$X_{q,\mu}/H_{k,u}:= \{ \{x\}*H_{k,u}:\> x\in X_{q,\mu}\} $$ (equipped with the
quotient topology)  carries a canonical quotient hypergroup structure 
with the convolution
$$ \delta_{ \{x\}*H_{k,u}} * \delta_{ \{y\}*H_{k,u}}= \int_{X_{q,\mu}} 
 \delta_{ \{z\}*H_{k,u}}\> d(\delta_x*\delta_y)(z) \quad\quad(x,y\in
 X_{q,\mu})$$
(where this convolution is independent of the representants $x,y$ of the
 cosets; this may fail for arbitrary subhypergroups of arbitrary commutative hypergroups).
Moreover, as in the group case, the hypergroup $X_{q,\mu}/H_{k,u}$
is isomorphic with $X_{q-k,\mu}$. This fact implies (see \cite{V2} and references
there)   that all
subhypergroups of $X_{q,\mu}$ have a number of nice analytic 
properties which are
obvious in the case of locally compact abelian  groups, but which 
 may  fail for
general commutative hypergroups. We therefore may say that our Bessel
 hypergroups on matrix cones are hypergroups  which is quite close 
to lca groups.
\end{remark}

In the end of this section, we classify all automorphisms for $\b F=\b R, \b C$.
For this we denote  the group of all hypergroup automorphisms of
 $X_{q,\mu}$
by $Aut(X_{q,\mu})$. For $\b F=\b R$ we prove:

\begin{theorem}\label{classiauto-reell}
 Let $\b F=\b R$ and $\mu>\rho-1$.
Then  $Aut(X_{q,\mu})=\{T_a:\> a\in GL(q)\}$. 
\end{theorem}

To describe all automorphisms for $\b F=\b C$ and $q\ge2$, $\mu>\rho-1$, we
denote the transposition $x\mapsto x^t$ on the space of Hermitian matrices by
$\tau$. We know from Eq.~(\ref{def-convo}) that its restriction to $\Pi_q$ is
contained in $Aut(X_{q,\mu})$. Moreover, for $q\ge2$, $\tau\not\in\{T_a:\>
a\in GL(q,\b C)\}$. (In fact, for the proof of this fact we may restrict our
attention to the case $q=2$, in which case the statement 
 can be checked by a direct computation.)

Moreover, as $\tau\circ T_a\circ \tau = T_{\tau(a)}$ for all $a\in GL(q,\b
C)$, it follows that 
$$\{\sigma\circ T_a:\>  a\in GL(q, \b C),\> \sigma\in\{Id,\tau\}\}$$
is a group of automorphisms of $X_{q,\mu}$.

\begin{theorem}\label{classiauto-complex}
 Let $\b F=\b C$ and $\mu>\rho-1$.
Then  
$$Aut(X_{q,\mu})=\{\sigma\circ T_a:\>  a\in GL(q, \b C),\> \sigma\in\{Id,\tau\}\}
.$$ 
\end{theorem}

Our proof of this classification is  quite complicated, covers
the remaining part of Section 4 and may be
skipped at a first reading. The proof
is divided into several steps which
partially work for all  $\b F$. We start with  $q=1$, i.e., 
 $X_{1,\mu}$ is a 
Bessel-Kingman hypergroup. This case was already handled in Zeuner \cite{Z1},
but we include the proof for sake of completeness.

\begin{lemma}\label{step1}
Let  $\mu$, $\b F$ arbitrary. Then $Aut(X_{1,\mu})=\{T_a:\> a > 0\}$.
\end{lemma}

\begin{proof} The convolution on $X_{1,\mu}=[0,\infty[$ satisfies 
$\{|a-b|, a+b\}\subset \{a\}*\{b\} \subset [|a-b|, a+b]$ for all $a,b\ge0$.
Therefore, if $T\in Aut(X_{1,\mu})$ satisfies $T(1)=c$ for some $c>0$,  we obtain
$T(1/n)=c/n$ for all $n\in \b Z_+$ and thus, $T(m/n)=cm/n$ for all $m,n\in \b
Z_+$.
Continuity then yields $T(x)=cx$ for all $x\ge0$.
\end{proof}

In the next steps we deal with the multi-dimensional case for arbitrary fields
$\b F$. The first main result will be Proposition \ref{step5} below where we show
that each $T\in Aut(X_{1,\mu})$, which preserves diagonal matrices, preserves
the norm for all matrices in $\Pi_q$.

\begin{lemma}\label{step2}
Let $\mu>\rho-1$ and  $T\in Aut(X_{q,\mu})$. Then for each $u\in U_q$ and
$k=0,\ldots,q$ there exists $\tilde u\in U_q$ with $T(H_{k,u})=H_{k,\tilde u}$.
\end{lemma}

\begin{proof}
Consider the maximal chain
$$\{0\}=H_{0,u}\subset H_{1,u}\subset \ldots \subset H_{q,u}= \Pi_q$$
of subhypergroups such that all inclusions are proper. The classification of all
subhypergroups
in Proposition \ref{all-subhyper} now leads to the claim.
\end{proof}

In the following, we denote the diagonal matrix  $diag(0,\ldots,
0,1,0,\ldots, 0)\in \Pi_q$ with $1$ at the $i$-th element by $e_i$
($i=1,\ldots,q$).

\begin{lemma}\label{step3}
Let $\mu>\rho-1$ and  $T\in Aut(X_{q,\mu})$. Then there exists $a\in Gl(q)$
such that $T_a\circ T(c\cdot e_{i})=c\cdot e_{i}$ for all $c\ge 0$
and $i=1,\ldots,q$. 
\end{lemma}

\begin{proof}
For $i=1,\ldots,q$ consider $r_i:=T(e_i)\in\Pi_q$. These matrices have rank 1
by Lemma \ref{step2}. We thus find vectors $x_i\in \b F^q$ with
$r_i^2=x_ix_i^*$. We claim that the $x_i$ are linearly independent.
In fact, if they would be dependent, we would find $x\in \b F^q\setminus
\{0\}$ with $r_ix=x_ix_i^*x=0$ for all $i$. In other words, $r_1,\ldots, r_q$
would be contained in a proper subhypergroup of $\Pi_q$. But this is
impossible by Lemma \ref{step2}, as $\Pi_q$ is the only subhypergroup
containing all $e_i$.
We thus see that the $x_i$ are linearly independent. Hence we find a unique
$a\in Gl(q)$ such that for all $i$,  $ax_i$ is the $i$th unit vector. 
This implies $T_a\circ T(e_{i})=e_{i}$ for all
$i$. We thus conclude from Lemma \ref{step2} that  $T_a\circ T$ is an
automorphism on the one-dimensional hypergroups $\{c\cdot e_{i}:\> c\ge0\}$
with $T_a\circ T(e_{i})=e_{i}$. Therefore, by Lemma \ref{step1},  $T_a\circ T$
is the identity on these subhypergroups. This proves the lemma.
\end{proof}

\begin{lemma}\label{step4}
Let   $T\in Aut(X_{q,\mu})$ with $T(c\cdot e_{i})=c\cdot e_{i}$ for all $c\ge 0$
and $i=1,\ldots,q$. Then $T(r)=r$ for all diagonal matrices $r\in\Pi_q$.
\end{lemma}

\begin{proof}
Let $T^*$ be the dual automorphism according to Proposition \ref{prop-automor}
We first fix $i =1,\ldots,q$, $c\ge0$ and $s\in\Pi_q$ and notice
$$\phi_{T^*(s)}(c\cdot e_{i})=\phi_s(T(c\cdot e_{i}))=\phi_s(c\cdot e_{i}).$$
The Taylor expansion (\ref{homexp}) of  $\mathcal J_\mu$ now yields
$$1-\frac{c^2}{4\mu} tr(e_i^2 T^*(s)^2)+O(c^4)= 
1-\frac{c^2}{4\mu}tr(e_i^2 s^2) +O(c^4)$$
for $c\to0$. As this holds for all $i$, the matrices $T^*(s)^2$
and $s^2$ have the same diagonal parts for any $s$.

 Now let $r=\sum_{i=1}^q
c_i e_i\in\Pi_q$ be an arbitrary diagonal matrix. Using 
$\phi_{T^*(cs)}(r)=\phi_{cs}(T(r))$ and  the  Taylor expansion (\ref{homexp}),
we obtain
$$ 1-\frac{c^2}{4\mu}tr(T(r)^2 s^2) +O(c^4) = J_\mu(\frac{1}{4}r^2T^*(cs)^2)$$
where, by our considerations above, $r^2T^*(cs)^2= c^2r^2s^2 + h$ for some
matrix $h=h(c,s,r)$ with zeros on the diagonal. Moreover, as $T^*(cs)$ as well
as $cs$ are both positive semidefinite, and as the absolute values of all
entries of a  positive semidefinite matrix are bounded by the maximum of the
absolute values of the diagonal entries, we have $h(c,r,s)=O(c^2)$ for
$c\to0$.
Hence,  again by (\ref{homexp}) ,
$$J_\mu(\frac{1}{4}r^2T^*(cs)^2)= J_\mu(\frac{1}{4}c^2r^2s^2) +O(c^4).$$
Combining all results, we obtain
$$1-\frac{c^2}{4\mu} tr(r^2 s^2)+O(c^4) =
J_\mu(\frac{1}{4}c^2r^2s^2) +O(c^4) = 
1-\frac{c^2}{4\mu}tr(T(r)^2 s^2)+O(c^4) $$
for $c\to0$. Hence,  $tr(r^2 s^2)=tr(T(r)^2 s^2)$ for all $s\in\Pi_q$ and thus
all  $s\in H_q$. 
As the trace forms a scalar product, we obtain $T(r)=r$ as claimed.
\end{proof}

\begin{proposition}\label{step5}
Let  $\mu>\rho-1$ and   $T\in Aut(X_{q,\mu})$. Then there exist $a\in Gl(q)$
and a mapping $h:U_q\to U_q$ with $h(I)=I$ such that 
$T_a\circ T(uru^*)=T_{h(u)}(uru^*)$
for all diagonal matrices $r\in\Pi_q$ and all $u\in U_q$. In particular, for
this $a$, we have $\|T_a\circ T(x)\|=\|x\|$ for all $x\in\Pi_q$.
\end{proposition}

\begin{proof} By Lemmas \ref{step3} and  \ref{step4} we find $a\in Gl(q)$
such that $T_a\circ T(r)=r$ for all diagonal matrices $r\in\Pi_q$. The same
argument together with change of basis show that for each $u\in U_q$ there
exists $h(u)\in GL(q)$ such that $T_a\circ T(uru^*)=T_{h(u)}(uru^*)$
for all diagonal matrices $r\in\Pi_q$. Taking $r=I$, we see that for any $u\in
U_q$,
$I=T_a\circ T(I)=T_a\circ T(uIu^*)=T_{h(u)}(uIu^*)=T_{h(u)}(I)$, and hence
$h(u)\in U_q$. The lemma is now obvious.
\end{proof}

We next restrict our attention to the case $q=2$ where we derive the
classification for $\b F=\b R,\b C$. The key will be:

\begin{lemma}\label{step6}
Let $q=2$ and  $\mu>\rho-1$. Then, for all $a,c\in ]0,\infty[$,
$$\left(\left\{ \left(\begin{array}{cc}  a &0\\0& 0\end{array}\right)\right\}
*\left\{ \left(\begin{array}{cc}  0 &0\\0& c \end{array}\right)\right\}\right)
\cap  \left\{r\in\Pi_q:\> \|r\|\ge \sqrt{a^2+c^2}\right\} 
= \left\{s(\beta):\>
|\beta|=1
\right\},$$
where, for $\beta\in\b F$ with $|\beta|=1$,
$$s(\beta):=\left(\begin{array}{cc}  a^2 & a\beta c \\ a\bar \beta c& c^2
  \end{array}\right)^{1/2}
=\frac{1}{a^2+c^2} 
\left(\begin{array}{cc}  a & \beta c \\ -\bar\beta c& a
  \end{array}\right)
\left(\begin{array}{cc}  \sqrt{a^2+c^2}  & 0\\ 0& 0 \end{array}\right)
\left(\begin{array}{cc}  a & -\beta c \\ \bar\beta c& a
  \end{array}\right).$$
\end{lemma}

\begin{proof} Let $v=\left(\begin{array}{cc}  \alpha & \beta  \\ \gamma & \delta
  \end{array}\right)\in M_2$ with $vv^*\le I$. Therefore, the $(1,1)$-entry 
$|\alpha|^2+|\beta|^2$ of $vv^*$ satisfies $|\alpha|^2+|\beta|^2\le 1$. Hence,
$|\beta|\le1 $.
Moreover, for each $\beta\in\b F$ with $|\beta|\le1 $ we obviously find a corresponding $v$
with $vv^*\le I$.
Using the convolution formula Eq.~(\ref{def-convo}) and the fact that
the matrix
$$\left(\begin{array}{cc}  a &0\\0& 0\end{array}\right)^2 + 
\left(\begin{array}{cc}  0 &0\\0& c \end{array}\right)^2 + 
\left(\begin{array}{cc}  a &0\\0& 0\end{array}\right)v
\left(\begin{array}{cc}  0 &0\\0& c \end{array}\right)+
\left(\begin{array}{cc}  0 &0\\0& c \end{array}\right)v^*
\left(\begin{array}{cc}  a &0\\0& 0\end{array}\right)
=\left(\begin{array}{cc}  a^2 & a\beta c \\ a\bar \beta c& c^2
  \end{array}\right)$$
has norm $(a^4+c^4+2a^2c^2|\beta|^2)^{1/2}$, the statements of the lemma now
follow easily.
\end{proof}

\begin{corollary}\label{step6a}
Let $T\in Aut(X_{2,\mu})$ with $T(r)=r$ for diagonal matrices $r\in
 \Pi_2$. Then, for all $a,c >0$ and $\beta\in\b F$ with $|\beta|=1$ and with
 the notion of Lemma \ref{step6}, $T(s(\beta))\in\{ s(\gamma):\> \gamma\in\b
 F, \> |\gamma|=1\}$.
 \end{corollary}

\begin{proof}
$T$ preserves the matrices 
$\left(\begin{array}{cc}  a &0\\0&0\end{array}\right)$ 
and $\left(\begin{array}{cc}  0 &0\\0&c\end{array}\right)$ for $a,c\ge0$. As
$T$ also preserves norms by Proposition \ref{step5}, the statement follows
from Lemma  \ref{step6}.
\end{proof}

Corollary \ref{step6a} now  leads easily to the claimed classification for $q=2$
and $\b F=\b R$.

\begin{proposition}\label{step7}
For  $\b F=\b R$,
 $Aut(X_{2,\mu})=\{T_a:\> a\in GL(2)\}$. 
\end{proposition}

\begin{proof} According to Lemmas \ref{step3} and \ref{step4}, it suffices to prove that 
any $T\in Aut(X_{2,\mu})$ with $T(r)=r$ for diagonal matrices $r\in\Pi_2$
has the form $T=T_u$ for some $u\in O(2)$.
For this take  $a,c\ge0$ and  $\beta=\pm 1$ and consider the matrices
$s(\beta)$ as above. Then, by Corollary  \ref{step6a}, 
 $T( s(\beta))\in \{s(1),s(-1)\}$.

Assume now that $s_0:=\left(\begin{array}{cc} 1 & 1 \\
 1&1\end{array}\right)$ satisfies $T(s_0)=s_0$.
The continuity of $T$ then implies that $T( s(\beta)) =s(\beta) $ for all
 $a,c$ and $\beta=1$. Clearly, this statement then must also hold for all $a,c$ and $\beta=-1$.
 As by  the diagonalization of $s(\beta)$ in Lemma \ref{step6}
each rank one matrix in $\Pi_2$ appears as some $s(\beta)$, we conclude 
 that $T$ is the identity for all rank one matrices.  Proposition \ref{step5} now
implies that  $T$ is the identity  for all matrices in $\Pi_q$.

Furthermore, if $T(s_0)= \left(\begin{array}{cc} 1 & -1 \\
- 1&1\end{array}\right)$, then we get by the same arguments $T( s(\beta)) =s(-\beta) $ for 
all  $a,c$ and $\beta$ and thus  $T=T_u$ for $u=\left(\begin{array}{cc} 1 & 0\\
0&-1\end{array}\right)$ on $\Pi_q$.
 This proves the claim.
\end{proof}

We next deal with $q=2$ and $\b F=\b C, \b H$. We here need the following
 observation:

\begin{lemma}\label{step7a}
Let $T\in Aut(X_{2,\mu})$ with $T(r)=r$ for all diagonal matrices
$r\in\Pi_2$. Then there exists $v\in U_2$ such that $T_v\circ T(r)=r$ for all
$r\in\Pi_2$ which are diagonal or which have the form 
$r= \left(\begin{array}{cc} s+t & s-t \\s-t&s+t\end{array}\right)$ with $s,t\ge0$.
\end{lemma}

\begin{proof} By Corollary \ref{step6a} there exist numbers $\beta_1(s),
  \beta_2(t)\in\b F$ with $|\beta_1(s)|=|\beta_2(t)|=1$ such that
\begin{equation}\label{st}
T\left(\left(\begin{array}{cc} s & s \\s&s\end{array}\right)\right)= 
\left(\begin{array}{cc} s & s\beta_1(s)
    \\s\overline{\beta_1(s)}&s\end{array}\right),\quad
T\left(\left(\begin{array}{cc} t & -t \\-t&t\end{array}\right)\right)= 
\left(\begin{array}{cc} t & -t\beta_2(t)
    \\-t\overline{\beta_2(t)}&t\end{array}\right)
\end{equation}
for all $s,t\ge0$. On the other hand, 
$$<\left(\begin{array}{cc} s & s \\s&s\end{array}\right),
\left(\begin{array}{cc} t & -t \\-t&t\end{array}\right)> \> =\> 
\{u_0ru_0^*:\> r\in\Pi_2 \>\>{\rm diagonal}\}=:H$$
for $u_0:=\frac{1}{\sqrt{2}}\left(\begin{array}{cc} 1 & 1
    \\1&-1\end{array}\right)\in U_2$, and thus, by Proposition
\ref{step5}, $T(w)=v^*wv$ for all $w\in H$ and some $v\in U_2$. In particular, $T$ is $\b
R$-linear on $H$ which ensures that $\beta_1,\beta_2$ in Eq.(\ref{st}) are constants
independent of $s,t\ge0$.
Therefore,
$$T\left(\left(\begin{array}{cc} s+t & s-t \\s-t&s+t\end{array}\right)\right)= 
\left(\begin{array}{cc} s+t & \beta_1s-\beta_2t \\\bar{\beta_1}s-\bar{\beta_2}t&s+t\end{array}\right)$$
for $s,t\ge0$. As $T$ is norm-preserving, it follows that
$|s-t|=|s-t\beta_1/\beta_2|$ for all $s,t\ge0$. This yields
$\beta_1=\beta_2=:\beta$, and thus $T(w)=v^*wv$ for 
$v=\left(\begin{array}{cc} \beta & 0\\0&\bar\beta\end{array}\right)\in U_2$,
Therefore, $v$ has the properties claimed in the lemma.
\end{proof}

We next derive the classification for $\b F=\b C$ and $q=2$.

\begin{proposition}\label{step7b}
For  $\b F=\b C$,
 $$Aut(X_{2,\mu})=\{\sigma\circ T_a:\> a\in GL(2,\b C), \> \sigma\in\{Id,\tau\}\} .$$
\end{proposition}

\begin{proof}
According to Lemmas \ref{step3} and \ref{step4} and \ref{step7a}, it suffices to consider
 $T\in Aut(X_{2,\mu})$ with $T(r)=r$ for all $r\in\Pi_2$ which are diagonal or
 which have the form 
$r= \left(\begin{array}{cc} s+t & s-t \\s-t&s+t\end{array}\right)$ with $s,t\ge0$.
Let $s,t\ge0$ and let $u_0:=\frac{1}{\sqrt{2}}\left(\begin{array}{cc} 1 & 1
    \\1&-1\end{array}\right)\in U_2$
with $u_0^{-1}=u_0$. As $T_{u_0}$ preserves the norm on $\Pi_2$, we conclude
 from Lemma 4.14 and the fact that $T_{u_0}$ is a hypergroup automorphism that
\begin{align}
\left(\left\{T_{u_0} \left( \left(\begin{array}{cc}  s &0\\0& 0\end{array}\right)\right)\right\}
*\left\{T_{u_0} \left( \left(\begin{array}{cc}  0 &0\\0&
        t\end{array}\right)\right)\right\}\right)
&\cap\{r\in\Pi_q:\> \|r\|\ge\sqrt{s^2+t^2}\}
\nonumber\\&=
\{t(\beta):\> \beta\in\b C,\> |\beta|=1\}
\end{align}
with
\begin{align}
t(\beta)&=T_{u_0} \left( \left(
\begin{array}{cc}  s^2 &st\beta\\st\bar\beta& t^2\end{array}\right)^{1/2}\right)
\nonumber\\&=\frac{1}{\sqrt{2}} \left(\begin{array}{cc} 
 s^2 +t^2+st(\beta+\bar\beta)&s^2 -t^2+st(\bar\beta-\beta)\\
s^2 -t^2+st(\bar\beta-\beta)&s^2 +t^2-st(\bar\beta+\beta)
 \end{array}\right)^{1/2}.\nonumber
\end{align}
As
$$T\circ T_{u_0} \left( \left(\begin{array}{cc}  s &0\\0&
      0\end{array}\right)\right)
=T_{u_0} \left( \left(\begin{array}{cc}  s &0\\0&
      0\end{array}\right)\right)
\>\>{\rm and}\>\>
T\circ T_{u_0} \left( \left(\begin{array}{cc}  0 &0\\0&
      t\end{array}\right)\right)
=T_{u_0} \left( \left(\begin{array}{cc}  0 &0\\0&
      t\end{array}\right)\right)$$
by our assumption, we see that for all $\beta\in\b C$ with $|\beta|=1$ there
      exists $\gamma=\gamma(\beta)\in\b C$ with $|\gamma|=1$ such that
      $T(t(\beta))=t(\gamma)$. On the other hand, $t(\beta)$ is a rank one
      matrix and has thus the form $s(\delta)$ for some $\delta\in\b C$ with
      $|\delta|=1$ in the notion of Lemma \ref{step6}. Therefore, by Corollary
 \ref{step6a}, the diagonal entries of $t(\beta)^2$ are preserved under
      $T$. Hence, $\beta$ and $\gamma(\beta)$ have the same real parts, and
      thus $T(t(\beta))\in\{t(\beta),t(\bar\beta)\}$. A continuity argument
      shows that we have either $T(t(\beta))=t(\beta)$ for all $s,t\ge0$ and
      all $|\beta|=1$, orthat we always have the other case. As each rank one matrix
      $r\in\Pi_2$ appears as some $t(\beta)$ (for suitable $s,t,\beta$), we conclude from Proposition
      \ref{step5} that $T$ is either the identity or the transposition $\tau$
      on $\Pi_2$.
\end{proof}

We next restate Propositions \ref{step7} and  \ref{step7b}. For this, 
we define for $i,j=1,\ldots,q$
 with $i\ne j$
 the space $U^{i,j}(q)$ of all unitary  $v\in U_q$ 
with $v_{k,k}=1$ for all $k\ne i,j$, i.e.,
 there are at most two  possible non-trivial
non-diagonal entries of $v$ in  the positions $(i,j)$ and $(j,i)$.
The following  statement now follows immediately from the proof of 
 Propositions \ref{step7} and  \ref{step7b} by a suitable basis change.

\begin{lemma}\label{step8}
Let $b F=\b R,\b C$,   $q\ge2$, $u\in U_q$, and $T\in Aut(X_{q,\mu})$ with
 $T(uru^*)=uru^*$ for all diagonal matrices $r$. Let  $i,j=1,\ldots,q$
 with $i\ne j$. Then 
 there exists  $u(i,j)\in U_q$  and $\phi\in\{Id, \tau\}$ such that 
$T(vuru^*v^*)=\phi\circ T_{u(i,j)}(vuru^*v^*)$ 
 for all  $v\in U^{i,j}(q)$ and all diagonal matrices $r\in \Pi_q$.
For $\b F=\b R$ only the case $\phi=Id$ appears.
\end{lemma}

We are now ready to complete the classification.

\begin{proof}[Proof of Theorems \ref{classiauto-reell} and
\ref{classiauto-complex}.]
 
It suffices to consider the slightly more complicated case $\b F=\b C$.
Moreover, as in the proof of Propositions \ref{step7} and  \ref{step7b}, it suffices to
prove that each $T\in Aut(X_{q,\mu})$ with $T(r)=r$ for
  diagonal matrices $r\in\Pi_q$ has the form  $T=T_u$ or $T=\tau\circ T_u$  for some $u\in U(q)$.
To prove this we recapitulate  that 
there exist $N=N(q)$ and $i_1,\ldots, i_N,j_1,\ldots, j_N\in\{1,\ldots, q\}$
with $i_n\ne j_n$ for $n=1,\ldots,N$ such that
$$U_q= U^{i_1,j_1}(q) \cdot U^{i_2,j_2}(q)\cdots U^{i_N,j_N}(q).$$
Moreover, Lemma \ref{step8} and induction show that for 
$n=0,1,\ldots,N$ there exist $u_n\in U_q$ and
$\phi_1,\ldots,\phi_N\in\{Id,\tau\}$  such that 
$$T(v_nv_{n-1}\ldots v_1 r v_1^*\ldots v_{n-1}^* v_n^* )=\phi_n\circ T_{u_n}
(v_nv_{n-1}\ldots v_1 r v_1^*\ldots v_{n-1}^* v_n^* )$$
for all diagonal matrices $r\in\Pi_q$ and all $v_1\in O^{i_1,j_1}(q),\ldots ,
v_n\in O^{i_n,j_n}(q)$. 
The theorem now follows for $n=N$.
\end{proof}

We finally consider the case $\b F=\b H$. By the preceding proof,
 the classification only depends on the computation of
$Aut(X_{2,\mu})$. We suggest that here the study of concrete additional
matrices as in Lemma \ref{step7a} and Proposition  \ref{step7b}  leads
to the following conjecture:

Consider the group $G_{\b H}$ of  automorphisms of
the field $\b H=<1,i,j,k>_{\b R}$ 
which fix the real line, which is generated by the 3 automorphisms
which switch two of the $i,j,k$ and change the sign of the third component.
In this $D_3$-case we then we have $|G _{\b H}|=24$, and we may let act   $G_{\b H}$  on $ X_{2,\mu}$ 
by using the same transformation in each component of a matrix. It can be
easily checked by  Eq.~(\ref{def-convo}) that  $G_{\b H}$ then forms a
group of hypergroup automorphisms on $ X_{2,\mu}$. Moreover, for $\tau\in
G_{\b H}$
and $a\in GL(q)$ we have $\tau^{-1}\circ T_a\circ
\tau=T_{\tau(a)}$. Therefore,  $\{\tau\circ T_a:\> a\in GL(q),\> \tau\in G_{\b
  H} \}$
forms a group of hypergroup automorphisms on $ X_{2,\mu}$.
We   expect that 
$$Aut(X_{q,\mu})=\{\tau\circ T_a:\> a\in GL(q),\> \tau\in G_{\b H} \}.$$

\section{Convolution semigroups and Wishart distributions}

In this section we first introduce  convolution
  semigroups and associated random walks on $ X_{q,\mu}$. This concept is
  well-known for commutative hypergroups; see \cite{BH}, \cite{ReV}, and
  references there. We shall see that in particular 
general so-called squared Wishart distributions form such convolution semigroups.
In this way several known results about Wishart distributions and Wishart processes
may be partially seen under a new light, see
  \cite{B},\cite{CL},\cite{Co},\cite{Di},\cite{GY}, \cite{FK},\cite{H}, and in
  particular,\cite{Ja} and \cite{Mu}. We  here notice that 
these Wishart   distributions will appear later as limits in a central limit
  theorem in the next section. 
We first recapitulate the notions of convolution
  semigroups and associated random walks.

\begin{definition}\label{defsemi} 
\begin{enumerate}\itemsep=-2pt
\item[\rm{(1)}]  A family $(\mu_t)_{t\ge 0} \subset M^1(X_{q,\mu})$ of
 probability measures on $ X_{q,\mu}$ is called
a (continuous) convolution semigroup on $X_{q,\mu}$, if $\mu_s*\mu_t=\mu_{s+t}$
 for all $s,t\ge 0$ with $\mu_0=\delta_e$,
 and if the mapping $[0,\infty[\to M^1(X_{q,\mu})$, $t\mapsto \mu_t$
 is weakly continuous.
\item[\rm{(2)}]  A convolution semigroup $(\mu_t)_{t\ge 0}$
is called Gaussian if
$$\lim_{t\to0} \frac{1}{t} \mu_t(K\setminus U)=0
\quad\quad\text{ for all open  subsets}\>\>
U\subset X_{q,\mu} \>\>{\rm with}\>\> 0\in U.$$
\item[\rm{(3)}] Let $(\mu_t)_{t\ge0}$ be a convolution semigroup on $ X_{q,\mu}$.
 A $ X_{q,\mu}$-valued time-homogeneous Markov process $(X_t)_{t\ge0}$
is called a Le\'vy process on $ X_{q,\mu}$ associated with $(\mu_t)_{t\ge 0}$,
if its transition probabilities satisfy
$$P(X_t\in A|\> X_s=x)= (\mu_{t-s}*\delta_x)(A)$$
for all $0\le s\le t$, $x\in  X_{q,\mu}$, and Borel sets $A\subset  X_{q,\mu}$.
By well-known  general principles for  Feller processes, 
  a Le\'vy process on $ X_{q,\mu}$ always admits a
version with rcll paths, and  a Levy process is Gaussian, i.e.,  is
associated with a Gaussian convolution semigroup, if and only if it admits a
version with continuous paths; see \cite{ReV}. 
\item[\rm{(4)}] Similar to the continuous case, we say that
 random walk $(S_n)_{n\ge0}$ on $X_{q,\mu}$ associated
  with a sequence $(\mu_n)_{n\ge1}\subset M^1(\Pi_q)$ is a Markov chain 
with initial distribution $P_{S_0}=\delta_0$ and the transition probabilities
\begin{equation}\label{hom}
P(S_n\in A|\> S_{n-1}=x)\> =\> (\delta_x*\mu_n)(A)
\end{equation}
for $n\ge 1$, $x\in\Pi_n$ and Borel sets $A\subset\Pi_q$.
If all $\mu_n$ are equal to some $\mu$, then $(S_n)_{n\ge0}$ is 
time-homogeneous, and we say that it is associated with the measure $\mu$.

It is easy to check  that for a random walk $(S_n)_{n\ge0}$ on $X_{q,\mu}$
associated  with $(\mu_n)_{n\ge1}$ 
  and   $n\ge0$,  $S_n$ has distribution
  $P_{S_n}=\mu_1*\mu_2*\cdots*\mu_n*P_{S_0}$.
\end{enumerate}
\end{definition}

We next turn to Wishart distributions which form examples of Gaussian
convolution semigroups. Before defining them, we point out at the beginning 
 that our notion of  
Wishart distributions is  equivalent to, but  slightly different from  the classical
one, as
here in the group  case $\mu=dp/2$, a positive semidefinite matrix
$r\in \Pi_q$ corresponds to $\sqrt{x^*x}$ for $x\in M_{p,q}$ and not 
to $x^*x$ as usual. 
In this way, images of $U_p$-invariant normal distributions on $M_{p,q}$
under the projection $M_{p,q}\to \Pi_q$, $x\mapsto \sqrt{x^*x}$, will be 
images of classical Wishart distributions on  $\Pi_q$ under  $r\mapsto \sqrt
r$ on  $\Pi_q$. 
Also for general
parameters
$\mu$, we use these images of classical Wishart distributions under
this square root mapping, and call these distributions 
squared Wishart distributions. For instance, for $q=d=1$,  classical
Wishart distributions are gamma distributions while 
squared ones are Rayleigh distributions. In this way, our
notion is in agreement with Kingman \cite{K} for $q=d=1$ and close to the classical
Euclidean setting. This notion has also the advantage that the limit theorems
in Sections 6 and 7 will be in our notion very close to the classical Euclidean setting.

\begin{definition}\label{standard-wishart}
The standard squared Wishart distribution $W=W(d,q,\mu)$ on $\Pi_q=\Pi_q(\b
F)$
with shape parameter $\mu\ge\rho-1$ is  the probability measure 
$${(2\pi)^{-q\mu}}e^{-tr(r^2)/2}\> d\omega_\mu(r) \quad\quad(r\in \Omega_q)$$
on $\Pi_q$.
This is fact a  probability measure; this follows for
instance from Lemma \ref{Fourier-general} below for $s=0$: 
\end{definition}

We next turn to general squared Wishart distributions and
observe first that  $\omega_\mu$ and hence  $W$  are invariant under the
unitary transforms $r\mapsto uru^*$ on $\Pi_q$ for $u\in U_q$.
As any $a\in M_q$ may be written as $a=su$ with $s=\sqrt{aa^*}\in \Pi_q$ and
$u\in U_q$,  the image $T_a(W)$ under the mapping 
$T_a(r)=\sqrt{ar^2a^*}$ agrees with $T_s(W)$, i.e., $T_a(W)$ depends
only on $s=\sqrt{aa^*}$.

\begin{definition}
The  squared Wishart distribution $W(s^2)=W(d,q,\mu;s^2)$ on $\Pi_q=\Pi_q(\b
F)$
with shape parameter $\mu\ge\rho-1$ and covariance $s^2$ for $s\in \Pi_q$ is
defined as image of $W$ under $T_s$ (or, by the preceding discussion, 
under $T_a$ for any $a\in M_q$ with $s^2=aa^*$).
\end{definition}

The transformation formula yields that for regular  $s\in \Pi_q$, the
 distribution  $W(s^2)$  has the $\omega_\mu$-density
\begin{equation}\label{density-Wishart}
f_{s^2}(r):= \frac{1}{\Delta(s)^{\mu} (2\pi)^{q\mu}}
e^{-tr(s^{-1}r^2s^{-1})/2}
\quad\quad (r\in\Pi_q).
\end{equation}
Moreover, if $s\in \Pi_q$ is singular with rank $k<q$, then $W(s^2)=T_s(W)$
is supported by  the proper subhypergroup $T_s(\Pi_q)$ which can be
identified with  $X_{k,\mu}$; cf. Proposition
\ref{all-subhyper}. We show below that if we regard 
$W(s^2)=T_s(W)$ as a measure on $X_{k,\mu}$, it again 
admits a density like Eq.~(\ref{density-Wishart}) with 
respect to the Haar measure
on $X_{k,\mu}$.

We next determine the  Fourier transforms of  squared Wishart distributions
(in the hypergroup sense).

\begin{lemma}\label{Fourier-general}
For $a\in M_q$, the  Fourier
 transform of  $W(aa^*)$
 is given by
$$\widehat{W(aa^*)}(s)= e^{-tr(a^*s^2a)/2} \quad\quad (s\in\Pi_q).$$
\end{lemma}

\begin{proof}
Proposition XV.2.1 of \cite{FK} yields 
$$\int_{\Omega_q} e^{-tr(xy)} {\mathcal J}_\mu(x) \Delta(x)^{\mu-n/q}\> dx =
\Gamma_{\Omega_q}(\mu)\Delta(y)^{-\mu} e^{-tr (y^{-1})}.$$
Change of variables $y^{-1}=s^2/2$ and $x=sr^2s/4$ the readily leads to the
claim for the standard case $a=I$. Finally Lemma \ref{Fourier-auto} shows that for any $a\in M_q$,
$$\widehat{W(aa^*)}(r)=\widehat{T_a(W)}(r)=\hat W(T_{a^*}(s))=
e^{-tr(a^*s^2a)/2}=e^{-tr(aa^*s^2)/2}.$$
\end{proof}

If we introduce the
exponentials
$e_z\in L^2(\Omega_\mu)$ with $e_z(r):=e^{-tr(r^2z^2)/2}$ for $a\in \Pi_q$, we may write the
preceding lemma briefly as 
\begin{equation}\label{e-squared}
 \widehat{e_z}=(2\pi)^{q\mu}\Delta(z)^{-\mu}e_{z^{-1}}.
\end{equation}

As announced above, we now briefly discuss the density of degenerated squared Wishart
distributions. We restrict our attention to a special case without loss of
generality (cf. Section \ref{quotient});
the general case would need to much additional notation.

%\begin{remark} 
Let $a=diag(1,\ldots,1,0,\ldots,0)$ be a diagonal matrix with
  rank $k\in\{0,\ldots,q\}$. Then the squared Wishart
distribution $W(a)$ is supported by  the subhypergroup $H_{k,I}$  which can be
identified with $X_{k,\mu}$ via 
$\left(\begin{array}{cc} r&0\\0&0\end{array}\right)\simeq r$; see
 Proposition  \ref{all-subhyper}. Moreover, characters of $H_{k,I}$
can be written as
$$r\in H_{k,I}\mapsto \phi_s(r)={\mathcal J}_\mu(\frac{1}{4}s^2r^2)=
{\mathcal J}_\mu\left(\frac{1}{4}
\left(\begin{array}{cc} s&0\\0&0\end{array}\right)^2
\left(\begin{array}{cc} r&0\\0&0\end{array}\right)^2\right)
\quad \quad(s\in H_{k,I});$$
see the proof of  Proposition  \ref{all-subhyper}. Therefore, 
for $s\in H_{k,I}$, 
$$\widehat{W(a)}(s)= \widehat{T_a(W)}(s)= \int_{H_{k,I}} \phi_s \> dT_a(W) = 
\int_{\Pi_q} \phi_s\circ T_a\> dW =
\widehat W\left(\begin{array}{cc} s&0\\0&0\end{array}\right).$$ 
Thus, by Lemma \ref{Fourier-general}, the Fourier transforms of $W(a)$ and the
standard squared Wishart distribution on $H_{k,I}$ are equal. The injectivity of the
hypergroup Fourier transform then yields that, under the identification above,
$W(a)$ is just the standard squared Wishart distribution on $H_{k,I}$.

We note that this result may be also obtained by direct computation.
The details here (e.g. regarding the constants) are however in our opinion more 
complicated than in our approach. Singular Wishart distributions are also
considered (after the transformation $r\mapsto\sqrt r$) in \cite{CL}, \cite{LM}.
%\end{remark}

We next collect some trivial properties of squared Wishart distributions.

\begin{lemma}\label{properties-wishart} For all $a,b\in\Pi_q$:
\begin{enumerate}
\item[\rm{(1)}] $T_a(W(b^2))=W(ab^2a)$;
\item[\rm{(2)}] $W(b^2)*W(a^2)=W(a^2+b^2)$;
\item[\rm{(3)}] $(\mu_t:=W(ta^2))_{t\ge0}$ is a Gaussian  convolution
  semigroup with $$\lim_{t\to 0} t^{-1} \mu_t(\Pi_q\setminus U)=0.$$
\end{enumerate}
\end{lemma}

\begin{proof} The first two statements follow from injectivity of the Fourier 
transform. Moreover, (2), the explicit formulas for the Fourier transforms,  and Levy's continuity theorem
for the hypergroup Fourier transform (Ch.~4.2 of \cite{BH}) imply that
$(\mu_t:=W(ta^2))_{t\ge0}$ is a convolution semigroup. For the proof  of being
Gaussian, we may assume $a$ as identity matrix, in which case the definition
may be checked  easily by using the transformation formula.
\end{proof}

We expect that
all Gaussian convolution semigroups on $X_{q,\mu}$ are given by squared
Wishart distributions in this way.
In fact, for the group cases $\mu=dp/2$ this  can be easily deduced from the
well-known corresponding result on the group $M_{q,p}\simeq \b R^{dpq}$.
In the other cases we shall investigate this point in a forthcoming paper.
We also point out that the generator of the Wishart semigroups is known; see 
Bru \cite{B}.

In the end of this section we  determine translates $\delta_x*W(s^2)$ of squared Wishart
distributions, as these shifted squared Wishart distributions appear in the
transition kernels of Gaussian  processes; cf. Section \ref{defsemi}(2). 
 We here follow ideas of C.~Herz \cite{H} and 
use the
following generalization of a result of Tricomi (see also \cite{Di}):

\begin{lemma}\label{tricomi}
Let $f,g\in L^2(\Omega_\mu)$, and $\phi(z):=\langle e_z,f\rangle$ and
 $\psi(z):=\langle e_z,g\rangle$ for $z\in\Omega_\mu)$. Then $\hat f=g$ if and
only if $\phi(z)=(2\pi)^{-q\mu}\Delta(z)^{-\mu}\psi(z^{-1})$.
\end{lemma}

\begin{proof}
Let $g=\hat f$.
Eq.~(\ref{e-squared}),
 and  the Plancherel
formula \ref{main2}(c) imply
\begin{align}
\phi(z):&=\langle e_z,f\rangle=(2\pi)^{-2q\mu}\langle\widehat{e_z},g\rangle
\nonumber\\&=(2\pi)^{-q\mu}\Delta(z)^{-\mu}\langle e_{z^{-1}}g\rangle
=(2\pi)^{-q\mu}\Delta(z)^{-\mu}\psi(z^{-1}).\nonumber
\end{align}
Conversely, a Stone-Weierstrass argument shows that the $e_z$ span a dense
subspace of $L^2(\Omega_\mu)$ which yields the converse statement.
\end{proof}

\begin{lemma}\label{translated-gaussian}
For any $s\in \Omega_q$ and $x\in\Pi_q$, 
$$d(\delta_x*W(s^2))(y)=\frac{1}{\Delta(s)^{\mu} (2\pi)^{q\mu}}
e^{-tr(x^2 +s^{-1}y^2s^{-1})/2} {\mathcal J}_\mu(-\frac{1}{4}x^2s^{-1}y^2s^{-1})\> d\omega_\mu(y)$$
\end{lemma}

\begin{proof}
Again, using a suitable automorphism, we may restrict our attention to the 
standard case. Fix $x\in\Pi_q$ and consider the functions
$f(s):=e^{-tr(s^2+x^2)/2}{\mathcal J}_\mu(-\frac{1}{4}x^2s^2)$ and
$g(r):= e^{-tr(r^2)/2}{\mathcal J}_\mu(\frac{1}{4}x^2r^2)$ in $L^2(\omega_\mu)$.
Then the associated functions $\phi,\psi$ according to the preceding lemma are
given by
$$\phi(z):=\int e^{-tr(s^2+x^2)/2}{\mathcal J}_\mu(-\frac{1}{4}x^2s^2)e^{-tr(s^2z^2)/2}\> d\omega(s)=
e^{-tr(x^2)/2} \widehat{e_{\sqrt{z^2+1}}}(ix)$$
and
$$\psi(z):=\int e^{-tr(r^2)/2} {\mathcal J}_\mu(-\frac{1}{4}x^2r^2) e^{-tr(r^2z^2)/2} \>
d\omega(r)=
 \widehat{e_{\sqrt{z^2+1}}}(x).$$
Eq.~(\ref{e-squared}) and analytic continuation yield
$\phi(z)=\Delta(z)^{-\mu}\psi(z^{-1})$, and hence, by Lemma \ref{tricomi},
$(2\pi)^{-q\mu}\hat f =g=(W*\delta_x)^\wedge$ on $\Pi_q$. 
As the hypergroup Fourier transform is injective, the proof is complete.
\end{proof}

\section{Limit theorems}

In this section we derive a central limit theorem as well as strong laws of
large numbers for random walks on 
matrix Bessel hypergroups which reduces in the group cases $\mu=pd/2$ just to radial
parts of the
classical central limit theorem on the vector space $M_{p,q}$ for sums 
of iid random variables and the corresponding classical strong laws of large
numbers of Kolmogorov. The proof of the central limit theorem
 is standard and uses a Taylor expansion of
the Fourier transforms as well as L\'evy's continuity theorem for hypergroups.
Before stating the CLT, we introduce so-called moment functions 
on matrix Bessel hypergroups. Such moment functions on hypergroups 
were introduced by Zeuner and used later for several limit theorems 
on hypergroups; see the monograph \cite{BH} for more details and references.

To introduce moment functions, we recapitulate that we 
regard $M_q$ and $H_q$ as real vector spaces  with 
 scalar product 
 $(x\vert y) := \mathfrak R {tr}(xy^*)$ and  norm $\|x\|=(x\vert x)^{1/2}$.
For $k\in \b Z_+$ and a function $g$ in the variable $s\in H_q$, we denote the
$k$-th differential of $g$ by $d_s^kg(s)$ where this is a $k$-linear map on
$H_q$.
 Following the literature on limit theorems on hypergroups (see Chapter 7 of
 \cite{BH}, \cite{Z2} and references cited there),
 we introduce moment functions and moments of probability measures  on
 $\Pi_m$.

\begin{definition}
For $k\in \b Z_+$  and $s_1,\ldots,s_k\in H_q$ define the moment function
$$m_k^{s_1,\ldots,s_k}(r):= i^k\cdot d_s^k \phi_s(r)|_{s=0}(s_1,\ldots,s_k)$$
on $\Pi_q$. As $\phi_s(r)={\mathcal J}_\mu(\frac{1}{4}sr^2s)$ and
${\mathcal J}_\mu(x)=1-\frac{tr(x)}{\mu}+o(\|x\|^2)$ for $x\to0$,
we have for $r\in \Pi_q$:
\begin{enumerate}\itemsep=-2pt
\item[\rm{(1)}] $m_{2k+1}^{s_1,\ldots,s_{2k+1}}=0$ for all $k$ and
  $s_1,\ldots,s_{2k+1}\in H_q$, and 
\item[\rm{(2)}] $m_2^{s_1,s_2}(r)= \frac{1}{2\mu} \mathfrak R tr(s_1r^2s_2).$
\item[\rm{(3)}]In particular,  $m_2^{I,I}(r)=\frac{1}{2\mu}\|r\|^2$.
\item[\rm{(3)}] $m_0\equiv 1$.
\end{enumerate}
 We say that a 
probability measure $\mu\in M^1(\Pi_q)$ has a
$k$-th moment, if $$\int_{\Pi_q} \|r\|^{k}\> d\mu(r)< \infty.$$
\end{definition}

In order to get estimates for moment functions and derivatives of $\hat\mu$,
we use the following Bochner-type integral representation for $\phi_s(r)$ 
 of \cite{R}:
\begin{equation}\label{Bochner-rep}
\phi_s(r) = \,\frac{1}{\kappa_\mu} \int_D e^{-i(rv\vert s)} \Delta(I-vv^*)^{\mu-\rho}dv.
\end{equation}

\begin{lemma}\label{moment-function-estimate}
For $k\in \b Z_+$  and $s_0,s_1,\ldots,s_k,r\in\Pi_q$,
$$\left| d_s^k \phi_s(r)|_{s=s_0}(s_1,\ldots,s_k)\right| \le
\|r\|^k\cdot\prod_{l=1}^k \|s_l\|,$$
and in particular, $|m_{2k}^{s_1,\ldots,s{2k}}|\le \|r\|^{2k}\cdot\prod_{l=1}^{2k} \|s_l\|$.
\end{lemma}

\begin{proof}
As $D$ is compact, we may interchange derivatives and integration in
Eq.~(\ref{Bochner-rep}).
Thus,
$$d_s^k \phi_s(r)|_{s=s_0}(s_1,\ldots,s_k) =\,\frac{(-i)^k}{\kappa_\mu}
 \int_D \prod_{l=1}^k (rv\vert s_l)\cdot e^{-i(rv\vert s)}
 \Delta(I-vv^*)^{\mu-\rho}dv.$$
As for $v\in D$, $0\le vv^*\le I$ and hence $0\le rvv^*r\le r^2$, we obtain
$$|(rv\vert s_l)|\le \|rv\|\cdot \|s_l\|\le \|r\|\cdot \|s_l\|,$$ and the lemma
 follows by taking absolute values.
\end{proof}

\begin{proposition}\label{moment-estimate}
Let  $\nu\in M^1(\Pi_q)$ and $k\ge1$. If the $2k$-th moment of $\nu$ exists,
then $\hat \nu$ is $2k$-times continuously differentiable on $\Pi_q$
with $d^{2k-1}\hat\nu(0)\equiv 0$.
Moreover, for $l\le 2k$,
$$d^l\hat \nu(s) = \int_{\Pi_q} d_s^l\phi_s(r)\> d\nu(r)
\quad\quad(s\in\Pi_q),$$
and, in particular for $l\le k$ and $s_1,\ldots,s_{2l}\in \Pi_q$,
$$d^{2l}\hat\nu(0)(s_1,\ldots,s_{2l})= \int_{\Pi_q}
m_{2l}^{s_1,\ldots,s_{2l}}(r)\> d\nu(r).$$
\end{proposition}

\begin{proof}
The estimation in Lemma 
\ref{moment-function-estimate} and standard results on derivatives of
parameter integrals ensure that partial differentiation up to order $2k$
and integration may be interchanged in $\hat\nu(s)=\int_{\Pi_q}\phi_s(r)\>
d\nu(r)$.
Therefore, under this condition, all statements are clear by the definition of
moment functions.
\end{proof}

We are now in the position to prove the following central limit theorem.
In the group case $\mu=dp/2$ it is equivalent to the classical central limit
theorem on the Euclidean space $M_{p,q}$ for sums of i.i.d. random variables
with a distribution which is invariant under the action of $U_p$ on
$M_{p,q}$.

\begin{theorem}\label{clt}
Let $\nu\in M^1(\Pi_q)$ such that the  $2k$-th moment of $\nu$ exists. 
Then, the matrix   $\sigma^2 :=\frac{1}{2\mu}  \int_{\Pi_q}
r^2\>d\nu(r)\in \Pi_q$
exists, and the probability measures 
$T_{n^{-1/2}I}(\nu^{(n)})=
(T_{n^{-1/2}I}(\nu))^{(n)}$ tend weakly to the squared Wishart distribution 
$W(\sigma^2)$ for $n\to\infty$.
\end{theorem}

\begin{proof} 
We first note that
\begin{align}
\int_{\Pi_q} m_2^{s,s}(r)\> d\nu(r)&= \frac{1}{2\mu} \int_{\Pi_q} tr(sr^2s)\>
d\nu(r)= \frac{1}{2\mu} tr\left(s\cdot  \int_{\Pi_q} r^2 d\nu(r)\cdot s\right)
\nonumber \\ &=tr(s\sigma^2s).\nonumber
\end{align}
Hence,  the preceding proposition and the  Taylor formula imply that $\hat\nu$ is twice continuously
  differentiable with 
$$ \hat\nu(s) = 1-\frac{1}{2} d^2\hat\nu(0)(s,s)+ o(\|s\|^2)
=
1-\frac{1}{2}tr(s\sigma^2s)+ o(\|s\|^2)$$
for $s\to0$. Hence,by Lemma \ref{Fourier-auto}, for any $s\in\Pi_q$,
$$\lim_{n\to\infty}
((T_{n^{-1/2}I}(\nu))^{(n)})^\wedge(s)=
\lim_{n\to\infty}\left( 1-\frac{1}{2n}tr(s\sigma^2s)+
  o(n^{-1})\right)^n=e^{-tr(s\sigma^2s)/2}.$$
Lemma \ref{Fourier-general} and Levy's continuity theorem for the hypergroup
Fourier transform (see Section 4.2 of \cite{BH}) now complete the proof.
\end{proof}

We next turn to strong laws of large numbers  for  random walks on
$X_{q,\mu}$. 
For this  we
use the algebraic properties of the moment functions on  $X_{q,\mu}$ which then
will be used to
construct martingales. 

\begin{lemma}
For all  $k\ge0$ and $s_1,\ldots,s_k,x,y\in\Pi_q$,
$$\int m_k^{ s_1,\ldots,s_k}\>
d(\delta_x*\delta_y)
= \sum_{l=0}^k \sum_{1\le i_1<i_2<\ldots<i_l\le k} m_l^{
  s_{i_1},\ldots,s_{i_l}}(x)
m_{k-l}^{s_1,\ldots,s_k\setminus s_{i_1},\ldots,s_{i_l}}(y)$$
where $s_1,\ldots,s_k\setminus s_{i_1},\ldots,s_{i_l}$ stands for the $s_j$
with $j\not\in \{i_1,\ldots,i_l\}$. In particular,
$$\int m_2^{ s_1,s_2}\>
d(\delta_x*\delta_y)
=m_2^{ s_1,s_2}(x)+m_2^{ s_1,s_2}(y)$$
and, in the language of matrix valued integrals,
$$\int r^2 \> d(\delta_x*\delta_y)(r)= x^2 +y^2.$$
\end{lemma}

\begin{proof}
As the  $\phi_s$ are multiplicative, the definition of moment functions yields
$$\int m_k^{ s_1,\ldots,s_k}\>
d(\delta_x*\delta_y) = i^k\cdot d_s^k
(\phi_s(x)\phi_s(y))|_{s=0}(s_1,\ldots,s_k).$$
The first  statement now follows readily from a multivariate version of the Leibniz
product rule for derivatives. The last equation follows from the preceding one
by taking matrices $s1,s_2\in\Pi_q$ which have only zero entries except for
precisely one 1 on the diagonal.
\end{proof}

We now construct martingales from random walks  $(S_n)_{n\ge0}$ on $X_{q,\mu}$ associated
  with a sequence $(\mu_n)_{n\ge1}\subset M^1(\Pi_q)$. For this we realize
  $(S_n)_{n\ge0}$
on some probability space which carries the canonical filtration associated
  with  $(S_n)_{n\ge0}$. Expectations will be denoted by $\b E$.

Based on the preceding observations, it is standard to derive the following
observations (cf. Section 7.3 of \cite{BH} or \cite{Z2}):

\begin{lemma}
\begin{enumerate}\itemsep=-2pt
\item[\rm{(1)}] For each $s\in \Pi_q$, the $\b R$-valued process 
 $$ \left(\phi_s(S_n)\cdot\prod_{k=1}^n (\hat\mu_k(s))^{-1}\right)_{n\ge0}$$ is a martingale.
\item[\rm{(2)}] Assume that all $\mu_n$ admit second moments. Then
for all $s_1,s_2\in\Pi_q$, 
$$\b E( m_2^{ s_1,s_2}(S_n))=
\sum_{k=1}^n \int m_2^{s_1,s_2}\> d\mu_k \quad\quad(n\ge0),$$
and 
$(m_2^{ s_1,s_2}(S_n)- \b E( m_2^{ s_1,s_2}(S_n)))_{n\ge0}$ is a martingale.

In matrix language, $\b E( S_n^2)=
\sum_{k=1}^n \int r^2 \> d\mu_k(r)$ for  $n\ge0$, and
$(S_n^2- \b E(S_n^2)_{n\ge0}$ is a matrix-valued martingale.
\end{enumerate}
\end{lemma}

Also higher moment functions can be used to construct martingales under
suitable moment conditions. This was worked out for instance in \cite{RV} for
the closely related case of Markov chains on Weyl chambers which are
associated with Dunkl operators. For a general discussion of moment functions
and associated martingales see also \cite{BH}, \cite{ReV}, \cite{Z2}.

Based on strong laws for martingales and the concept of moment functions, 
Zeuner \cite{Z2} derived general strong
laws of large numbers for random walks on general commutative hypergroups; see
also Section 7.3 of \cite{BH}. In the present setting, Zeuner's results
lead to the following strong laws which correspond in the group case
$\mu=pd/2$ precisely to the classical strong laws of large numbers of
Kolmogorov on the Euclidean spaces $M_{m,p}$.

\begin{theorem}
Let $(S_n)_{n\ge0}$ be a random walk on $X_{q,\mu}$ associated with the measusures
$(\mu_n)_{n\ge1}\subset M^1(X_{q,\mu})$.
\begin{enumerate}\itemsep=-2pt
\item[\rm{(1)}]  If $(a_n)_{n\ge1}\subset]0,\infty[$ satisfies $a_n\to\infty$
  and
$$\sum_{n=1}^\infty \frac{1}{a_n^2}E\left(\int\|r\|^2\> d\mu_n(r)\right)\> <\>
\infty,$$
then $\lim_{n\to\infty} S_n/a_n=0$ almost surely.
\item[\rm{(2)}] Let $\lambda\in]0,2[$ and  $(S_n)_{n\ge0}$  time-homogeneous
 with $\int \|r\|^2\> d\mu(r)<\infty$. Then $n^{-1/\lambda}S_n\to 0$ for $n\to\infty$
almost surely.
\end{enumerate}
\end{theorem}

\begin{proof}
Apply Theorem 7.8 and Corollary 7.11 of \cite{Z2} respectively  to the moment function 
$m_2^{I,I}(r)=\frac{1}{2\mu}\|r\|^2$.
\end{proof}

\begin{remark}
Let $\mu\ge \rho-1$. The compact group $U_q\subset GL(q)$ acts as group of
automorphisms on the hypergroup $X_{q,\mu}$
such that the space $X_{q,\mu}^{q_m}$ of orbits may be identified with the
Weyl chamber
$$W_q:=\{(\lambda_1,\ldots,\lambda_q)\in \b R^q:\>0\le\lambda_1\le\lambda_2\le\ldots\le\lambda_q\}$$
of type $B_q$. It is shown in  \cite{R} that $X_{q,\mu}^{U_q}\simeq W_q$
carries a commutative orbit hypergroup structure whose characters are
symmetric Dunkl kernels of type $B_q$; cf. \cite{Du}.  With the canonical
projection from $X_{q,\mu}$ onto $X_{q,\mu}^{U_q}\simeq W_q$, the preceding
limit theorems can be immediately be transferred into limit theorems for
random walks on this $B_q$-Dunkl-type hypergroup structure on $W_q$.
This leads to connections with limit results in \cite{RV}.
\end{remark}

\end{document}